\documentclass{elsart}
\usepackage{graphicx,epsfig}
\usepackage[cmtip,arrow]{xy}
\usepackage{pb-diagram,pb-xy}
\usepackage{subfigure,amsmath,amssymb}
\journal{Physica D}

\newcommand{\pfrac}[2]{\frac{\partial #1}{\partial #2}}
\newcommand{\ppfrac}[3]{\frac{\partial^2 #1}{\partial #2 \partial #3}}

\newcommand{\undemi}{\frac{1}{2}}
\newenvironment{disarray}
 {\everymath{\displaystyle\everymath{}}\array}
 {\endarray}
\newcommand{\qz}{q_0}
\newcommand{\qn}{q_n}
\newcommand{\qk}{q_k}
\newcommand{\qko}{q_{k+1}}
\newcommand{\qkmo}{q_{k-1}}

\newcommand{\qdk}{q^d_k}

\newcommand{\pz}{p_0}
\newcommand{\pn}{p_n}
\newcommand{\pk}{p_k}
\newcommand{\pko}{p_{k+1}}

\newcommand{\pdk}{p^d_k}
\newcommand{\tpdk}{\tilde{p}^d_k}

\newcommand{\Qdk}{Q^d_k}

\newcommand{\Pdk}{P^d_k}
\newcommand{\tPdk}{\tilde{P}^d_k}

\newcommand{\zk}{z_k}
\newcommand{\zko}{z_{k+1}}

\newcommand{\zdk}{z^d_k}

\newcommand{\Qk}{Q_k}
\newcommand{\Qko}{Q_{k+1}}

\newcommand{\Pk}{P_k}
\newcommand{\Pko}{P_{k+1}}

\newcommand{\Dt}{\Delta_\tau}
\newcommand{\Dtd}{\Delta_\tau^d}

\newcommand{\M}{{\mathcal M}}
\newcommand{\fdk}{f_k^d}
\newcommand{\gdk}{g_k^d}
\newcommand{\udk}{u_k^d}
\newcommand{\xdk}{x_k^d}
\newcommand{\xk}{x_k}
\newcommand{\xko}{x_{k+1}}
\newcommand{\hq}{ \Delta q}
\newcommand{\hp}{\Delta p}

\newcommand{\tk}{t_k}
\newcommand{\tko}{t_{k+1}}
\newcommand{\tkmo}{t_{k-1}}
\newcommand{\tdk}{t^d_k}
\newcommand{\edk}{e^d_k}
\newcommand{\ek}{e_k}
\newcommand{\eko}{e_{k+1}}

\newcommand{\cdk}{c^d_k}
\newcommand{\ck}{c_k}

\newcommand{\Cdk}{C^d_k}

\newcommand{\yk}{y_k}
\newcommand{\uk}{u_k}
\newcommand{\uko}{u_{k+1}}
\newcommand{\uok}{u_{1,k}}
\newcommand{\utk}{u_{2,k}}

\begin{document}

\begin{frontmatter}
\title{Discrete variational principles and Hamilton-Jacobi theory
for mechanical systems and optimal control problems.\thanksref{title}}
\thanks[title]{Research partially supported by NSF grants DMS $0103895$ and $0305837$.}
\author[aero]{V.M. Guibout\corauthref{info}},
\corauth[info]{Corresponding author.}
\ead{guibout@umich.edu}
\author[math]{A. Bloch}
\ead{abloch@umich.edu.}

\address[aero]{Department of Aerospace Engineering, University of Michigan, $1320$ Beal Avenue, Ann Arbor, MI 48109-2140}
\address[math]{Department of Mathematics, University of Michigan, $2074$ East Hall, Ann Arbor, MI 48109-1109}
\begin{abstract}
In this paper we present a general framework that allows one to study discretization of certain
dynamical systems. This generalizes earlier work on discretization of Lagrangian and Hamiltonian
systems on tangent bundles and cotangent bundles respectively. In particular we show how to obtain
a large class of discrete algorithms using this geometric approach. We give new geometric insight
into the Newmark model for example and  we give a direct discrete formulation of the
Hamilton-Jacobi method. Moreover we extend these ideas to deriving a discrete version of the
maximum principle for smooth optimal control problems.

We define discrete variational principles that are the discrete counterpart of known variational
principles. For dynamical systems, we introduce principles of critical action on both the tangent
bundle and the cotangent bundle. These two principles are equivalent and allow one to recover most
of the classical symplectic algorithms. In addition, we prove that by increasing the dimensionality
of the dynamical system (with time playing the role of a generalized coordinate), we are able to
add conservation of energy to any (symplectic) algorithms derived within this framework. We also
identify a class of coordinate transformations that leave the variational principles presented in
this paper invariant and develop a discrete Hamilton-Jacobi theory. This theory allows us to show
that the energy error in the (symplectic) integration of a dynamical system is invariant under
discrete canonical transformations. Finally, for optimal control problems we develop a discrete
maximum principle that yields discrete necessary conditions for optimality. These conditions are in
agreement with the usual conditions obtained from Pontryagin maximum principle. We illustrate our
approach with an example of a sub-Riemannian optimal control problem as well as simulations that
motivate the use of symplectic integrators to compute the generating functions for the phase flow
canonical transformation.

\end{abstract}

\end{frontmatter}

\section{Introduction}

Standard methods (called numerical integrators) for simulating motion take an initial condition and
move objects in the direction specified by the differential equations. These methods do not directly satisfy the physical conservation laws associated with the system. An alternative approach to
integration, the theory of geometric integrators\cite{mcl98qui,bud99ise}, has been developed over
the last two decades. These integrators strictly obey some of these physical laws, and take their
name from the law they preserve. For instance, the class of energy-momentum integrators conserves
energy and momenta associated with ignorable coordinates. Another class of geometric integrators is
the class of symplectic integrators which preserves the symplectic structure. This last class is of
particular interest when studying Hamiltonian and Lagrangian systems since the symplectic structure
plays a crucial role in these systems\cite{blo03bail,abr78,arn88a}. The work done by
Wisdom\cite{wis91,wis92} on the $n$-body problem perfectly illustrates the benefits of such
integrators.

At first, symplectic integrators were derived mostly as a subclass of Runge-Kutta algorithms for which the Runge-Kutta coefficients satisfy specific relationships \cite{san94cal}. Such a methodology, though very systematic, does not provide much physical insight and may be limited when we require several laws to be conserved. Other methods were developed in the $90$'s, among which we may cite the use of generating functions for the canonical transformation induced by the phase flow\cite{cha90sco,ge95wan} and the use of discrete variational principles. This last method ``gives a comprehensive and unified view on much of the literature on both discrete mechanics as well as integration methods''(Marsden and West\cite{mar01}). Names of variational principles differ in the literature, so we have decided to refer to Goldstein\cite{gol80} in this paper: Hamilton's principle concerns Lagrangian systems (i.e., refers to a principle of critical action that involves the Lagrangian) whereas the modified Hamilton's principle concerns Hamiltonian systems (i.e., refers to a principle of critical action that involves the Hamiltonian). Several versions of the discrete modified Hamilton's principle can be found in the literature such as the one developed by Shibberu\cite{shi92} and Wu\cite{wu89}. For the discrete Hamilton's principle, Moser and Veselov\cite{mos91ves} and then Marsden, West and Wendlandt \cite{mar01,wen97} developed a fruitful approach. Also, Jalnapurkar, Pekarsky and West \cite{jal00pek} developed a variational principle on the cotangent bundle based on generating function theory.

 In this paper, we focus on the discrete variational principles introduced by Guo, Li and Wu \cite{guo02li1,guo02li2,guo02li3} because the theory they have developed provides both a discrete modified Hamilton's principle (DMHP) and a discrete Hamilton's principle (DHP) that are equivalent. We modify and generalize both variational principles they introduce by changing the time discretization so that a suitable analogue of the continuous boundary conditions may be enforced. These boundary conditions are crucial for the analysis of optimal control problems and play a fundamental role in dynamics. Our approach not only allows us to obtain a large class of discrete algorithms but it also gives new geometric insight into the Newmark model \cite{new59}. Most importantly, using our improved version of the discrete variational principles introduced by Guo et al., we develop a discrete Hamilton-Jacobi theory that yields new results on symplectic integrators.

In the first part of this paper (sections \ref{sec:dcriticalactions}, \ref{sec:dcomparison} and
\ref{sec:examples}), we present a discrete Hamilton's principle on the tangent bundle and a discrete
modified Hamilton's principle on the cotangent bundle (section \ref{sec:dcriticalactions}), we
discuss the differences with other works on variational integrators (section \ref{sec:dcomparison})
and show that we are able to recover classical symplectic schemes (section \ref{sec:examples}). The
second part (sections \ref{sec:energy} and \ref{sec:dhj}) is devoted to issues related to energy
conservation and energy error. We first show that by considering time as a generalized coordinate
we can ensure energy conservation (section \ref{sec:energy}). Then we introduce the framework for
discrete symplectic geometry and the notion of discrete canonical transformations. We obtain a
discrete Hamilton-Jacobi theory that allows us to show that the energy error in the symplectic
integration of a dynamical system is invariant under discrete canonical transformations (section
\ref{sec:dhj}). Finally, in the last part (section \ref{sec:optimal}) we develop
a discrete maximum principle that yields discrete necessary conditions for optimality. These
conditions are in agreement with the usual conditions obtained from Pontryagin maximum principle
and define symplectic algorithms that solve the optimal control problem.


In each part, we illustrate some of the ideas with simulations. In particular we show in the first
part that symplectic methods allow one to recover the generating function from the phase flow while
standard numerical integrators fail because they do not enforce the necessary exactness condition.
The examples presented are the simple harmonic oscillator and a nonintegrable system describing a
particle orbiting an oblate body. In the second part we look at the energy error in the integration
of the equations of motion of a particle in a double well potential using a set of coordinates and
its transform under discrete symplectic map. In the last part, we use the discrete maximum
principle to study the Heisenberg optimal control problem.

\section{Discrete principles of critical action: DMHP and DHP}
\label{sec:dcriticalactions}
In this section, we develop a modified version of both variational principles introduced by Guo, Li
and Wu \cite{guo02li1,guo02li2,guo02li3} and present the geometry associated with them.

\subsection{Discrete geometry}
Consider a discretization of the time $t$ into $n$ instants ${\mathcal T}=\{(t_k)_{k\in [1,n]}\}$.
Here $\tko-\tk$ may not be equal to $\tk-\tkmo$ but for sake of simplicity we assume in the
following that $\tko-\tk=\tau$ $\forall k\in[1,n]$. The configuration space at $t_k$, is the
$n$-dimensional manifold $M_k$ and ${\mathcal M}=\bigcup M_k$ is the configuration space on
${\mathcal T}$. Define a discrete time derivative operator $\Dtd$ on ${\mathcal T}$.
Note that $\Dtd$ may not verify the usual Leibnitz law but a modified one. For instance, if we
choose $\Dtd$ to be the forward difference operator on $T{\mathcal T}$:
$$\Dtd q(t_k):=\frac{1}{\tau}(q(t_k+\tau)-q(t_k))=\frac{\qko-\qk}{\tau}:=\Dt q_k$$
then $\Dtd$ verifies:
\begin{equation}
\Dtd (f(t) g(t))=\Dtd f(t)\cdot g(t)+f(t+\tau)\cdot\Dtd g(t)\,.
\label{eq:modified_leibnitz}
\end{equation}

\subsection{Discrete Hamilton's principle}
Our modified version of the discrete Hamilton's principle derived by Guo, Li and Wu \cite{guo02li1}
is the discrete time counterpart of Hamilton's principle for Lagrangian systems. Consider a
discrete curve of points $(\qk)_{k\in [0,n]}$ and a discrete Lagrangian $L_d(\qdk,\Dtd \qdk)$ where
$\Dtd$ is a discrete time derivative operator and $\qdk$ is a function of $(\qk,\qko)$.

\begin{defn}[Discrete Hamilton's principle]
\label{defi:discrete_ham_pcp}
Trajectories of the discrete Lagrangian system $L_d$ going from $(t_0,q_0)$ to $(t_n,q_n)$ correspond to critical points of the discrete action
\begin{equation}
S^L_d=\sum_{k=0}^{n-1} L_d(\qdk,\Dtd \qk)\tau\,,
\end{equation}
in the class of discrete curves $(\qdk)_k$ whose ends are $(t_0,q_0)$ and
$(t_n,q_n)$. In other words, if we require that the variations of the discrete action $S^L_d$ be zero for
any choice of $\delta \qdk$, and $\delta q_0=\delta q_n=0$, then we obtain discrete Euler-Lagrange
equations.
\end{defn}
Note that if we do not impose $\tko-\tk=\tk-\tkmo$, then the discrete action would be defined as:
\begin{equation}
S^L_d=\sum_{k=0}^{n-1} L_d(\qdk,\Dtd \qk)(\tko-\tk)\,,
\end{equation}
but the discrete Hamilton's principle would be stated in the same manner\footnote{In this formulation,
the $t_k$'s are known, so there are no additional variables.}.

To proceed to the derivation of the equations of motion, we need to specify the derivative
operator, $\Dtd$. As we will explain below, its definition depends on the scheme we consider. We
should also mention that our variational principle differs from Guo, Li and Wu's since we consider
that the action has only finitely many terms and we impose fixed end points. Such a formulation is
more in agreement with continuous time variational principles and preserves the fundamental role
played by boundary conditions. For a discussion on this topic, we refer to Lanczos \cite{lan77}
section $15$.

\subsection{Discrete modified Hamilton's principle}
As in the continuous case, there exists a discrete variational principle on the cotangent bundle
that is equivalent to the above discrete Hamilton's principle.
\begin{defn}
Let $L_d$ be a discrete Lagrangian on $T{\mathcal M}$ and define the discrete Legendre transform
(or discrete fiber derivative) $\mathbb{F}L:T{\mathcal M}\rightarrow T^*{\mathcal M}$ which maps
the discrete state space $T{\mathcal M}$ to $T^*{\mathcal M}$ by
\begin{equation}
(\qdk,\Dtd \qdk)\mapsto(\qdk,\pdk)\,,
\end{equation}
where
\begin{equation}
\pdk=\pfrac{L_d(\qdk,\Dtd \qdk)}{\Dtd \qdk}\,.
\label{eq:legendre}
\end{equation}
If the discrete fiber derivative is a local isomorphism, $L_d$ is called regular and if it is a
global isomorphism we say that $L_d$ is hyperregular.
\end{defn}
If $L_d$ is hyperregular, we define the corresponding discrete Hamiltonian function on
$T^*{\mathcal M}$ by
\begin{equation}
H_d(\qdk,\pdk)=\langle\pdk, \Dtd \qdk\rangle - L_d(\qdk,\Dtd\qdk)\,,
\label{eq:discrete_ham}
\end{equation}
where $\Dtd \qdk$ is defined implicitly as a function of $(\qdk,\pdk)$ through equation
(\ref{eq:legendre}). Let $S_d^H$ be the discrete action summation:
\begin{equation}
S^H_d=\sum_{k=0}^{n-1} \left(\langle\pdk, \Dtd\qdk \rangle-H_d(\qdk,\pdk)\right)\tau\,,
\end{equation}
where $\tau$ is to be replaced by $\tko-\tk$ if $\tko-\tk\ne \tk-\tkmo$. Then the discrete
principle of least action may be stated as follows:
\begin{defn}[Discrete modified Hamilton's principle]
Trajectories of the discrete Hamiltonian system $H_d$ going from $(t_0,q_0)$ to $(t_n,q_n)$ correspond to critical points of the discrete action $S_d^H$ in the class of discrete curves $(\qdk,\pdk)$
whose ends are $(t_0, q_0)$ and $(t_n,q_n)$.
\label{def:discrete_pcp_la}
\end{defn}
Again, for deriving the equations of motion we need to specify the discrete
derivative operator, $\Dtd$ and its associated Leibnitz law. It will
generally depend upon the
scheme we consider as we will see through examples later.

\section{Comparison with other classical variational principle}
\label{sec:dcomparison}
At this point it is of interest to compare discrete variational principles introduced in this paper
and other classical discrete variational principles. As we mentioned above, the discrete
variational principles we develop are inspired by the work of Guo, Li and Wu \cite{guo02li1} and we
explained above the key difference between our work and this earlier work. We now point out the
main differences of the work discussed here with that of Marsden and West, based on the variational
principle introduced by Moser and Veselov. In the following, DVPI refers to the discrete
variational principle developed by Moser, Veselov, Marsden, Wendlandt et al. whereas DVPII denotes
the discrete variational principles developed by Guo and this paper.

The first main difference lies in the geometry of both variational principles. Whereas the discrete
Lagrangian is a functional on $Q\times Q$ where $Q$ is the configuration space in DVPI, it is a
functional on $TQ$ in DVPII. As a consequence, DVPII has a form more like
that of the continuous case but has
a major drawback: we have to specify the derivative operator and the Leibnitz law it verifies in
order to derive discrete Euler-Lagrange equation. Such a law allows us to perform the discrete
counterpart of the integration by parts and depends on the scheme we consider. On the other hand,
the Euler-Lagrange equation obtained by DVPI is scheme independent. One benefit is that these
equations ensure satisfaction of physical laws such as Noether's theorem for any numerical scheme
which can be derived from them.

The next important difference between the two discrete variational principles lies in the role of
the Legendre transformation in defining a discrete Hamiltonian function from the discrete
Lagrangian. In DVPI, one defines a discrete Legendre transform to compute the momenta from the
discrete Lagrangian function, so one may study the discrete dynamics on both $Q\times Q$ and
$T^*Q$. However, it does not seem possible to define a discrete Hamiltonian function from the
discrete Lagrangian and develop a DMHP. Given a Hamiltonian system, to derive discrete equations
of motion using DVPI one needs to first find a continuous Lagrangian function by performing a
Legendre transform on the continuous Hamiltonian function, then apply DVPI and finally use the
discrete Legendre transform to study the dynamics on $T^*Q$ (see for instance
\cite{mar01} page $408$). While this point may not be of importance when dealing with dynamical systems, it is crucial
if one wants to discretize an optimal control problem, where the continuous Hamiltonian function
does not have any physical meaning and the Legendre transformation may not be well-defined (See
section \ref{sec:optimal}). DVPII naturally defines a discrete Legendre transform and a DMHP.

As mentioned in the introduction, people have already introduced DMHPs on the cotangent bundle,
but, as far as we know, no one has developed an approach that allows one to equivalently consider both the
Hamiltonian and Lagrangian approaches in discrete settings (i.e., a DMHP and a
DHP that are equivalent for non-degenerate Lagrangian systems). In addition, the DMHPs that can be found in the literature do not allow one
to recover most of the classical schemes. For instance, Shibberu's DMHP focuses on the midpoint
scheme and Wu developed a different DMHP for each scheme.

Let us now look at some classical schemes and see how they can be derived from DVPII.

\section{Examples}
\label{sec:examples}
\subsection{St\"{o}rmer's rule and Newmark methods}
St\"{o}rmer's scheme is a symplectic algorithm that was first derived for molecular dynamics problems.
It can be viewed as a Runge-Kutta-Nystr{\"o}m method induced by the leap-frog partitioned
Runge-Kutta method\cite{san94cal}. The derivation of St\"{o}rmer rule as a variational integrator came
later and can be found in \cite{wu89,wen97}. Guo, Li and Wu \cite{guo02li3} recovered this
algorithm using their discrete variational principles. In the next subsection, we briefly go
through the derivation and add to their work the velocity Verlet \cite{swo82} and Newmark
methods\cite{mar01}. In particular, we will show how the conservation of the Lagrangian and
symplectic two-form is built into DVPII.

\subsubsection{From the Lagrangian point of view}
We first let $\qdk=\qk$ and define the discrete Lagrangian by $L_d(\qdk,\Dtd\qk)=L(\qk,\Dtd\qk)$
and the discrete derivative operator as the forward difference $\Dtd =\Dt$. $\Dt$ satisfies the
modified Leibnitz law (\ref{eq:modified_leibnitz}).
Discrete equations of motion are obtained from discrete Hamilton's principle (definition
(\ref{defi:discrete_ham_pcp})):

\begin{eqnarray}
\delta S^L_d&=&\tau\sum_{k=0}^{n-1} \delta L_d(\qk,\Dt \qk)\\
&=&\tau\sum_{k=0}^{n-1} \langle D_1 L_d(\qk,\Dt \qk),\delta \qk\rangle +
\langle D_2 L_d(\qk,\Dt \qk),\delta \Dt \qk\rangle\label{eq:del_trans0}\\
&=&\tau\sum_{k=1}^{n-1} \langle D_1 L_d(\qk,\Dt \qk)-\Dt D_2 L_d(\qkmo,\Dt \qkmo),\delta \qk\rangle\nonumber\\
&&+\Dt \langle D_2 L_d(\qkmo,\Dt \qkmo),\delta \qk\rangle\nonumber\\
&&+\tau \langle D_1L_d(\qz,\Dt\qz) \delta \qz\rangle +\tau D_2 L_d(\qz,\Dt\qz)\delta \Dt \qz \label{eq:del_trans1}\\
&=&\tau\sum_{k=1}^{n-1} \langle D_1 L_d(\qk,\Dt \qk)-\Dt D_2 L_d(\qkmo,\Dt \qkmo),\delta \qk\rangle - \nonumber\\
&& - \langle D_2 L_d(\qz,\Dt\qz),\delta\qz \rangle+\tau \langle D_1L_d(\qz,\Dt\qz) \delta
\qz\rangle\nonumber\\&&+\langle D_2 L_d(q_{n-1},\Dt q_{n-1}),\delta
\qn\rangle\label{eq:del}\,,
\end{eqnarray}
where the commutativity of $\delta$ and $\Dt$ and the modified Leibnitz law defined by equation
(\ref{eq:modified_leibnitz}) have been used.

Discrete Euler-Lagrange equations follow by requiring the variations of the action to be zero for
any choice of $\delta \qk$, $k\in[1,n-1]$ and $\delta \qz=\delta \qn=0$:
\begin{equation}
D_1 L_d(\qk,\Dt \qk)-\Dt D_2 L_d(\qkmo,\Dt \qkmo)=0\,.
\label{eq:stormer_del}
\end{equation}
Suppose $L(q,\dot{q})=\undemi \dot{q}M\dot{q}-V(q)$, then equation (\ref{eq:stormer_del}) yields
St\"{o}rmer's rule:
\begin{equation}
\qko=2\qk-\qkmo+h^2 M^{-1}(-\nabla V(\qk))\,.
\label{eq:stormer_l}
\end{equation}

Consider the one-form\footnote{Einstein's summation convention is assumed}
$$\theta_k^L=\pfrac{L_d(\qkmo,\Dt\qkmo)}{\Dt\qkmo^i}d\qk^i\,,$$ and define
the Lagrangian two-form $\omega_k^L$ on $T_{q_k}\M$:
\begin{eqnarray}
\omega_k^L&=&d\theta_k^L\nonumber\\
&=&\ppfrac{L_d(\qkmo,\Dt \qkmo)}{\qkmo^i}{\Dt \qkmo^j} d\qk^i \wedge
d\qk^j+\ppfrac{L_d(\qkmo,\Dt\qkmo)}{\Dt\qkmo^i}{\Dt \qkmo^j} d\Dt\qk^i \wedge d\qk^j\,.
\end{eqnarray}

\begin{lem}
The algorithm defined by St{\"o}rmer's rule preserves the Lagrangian two-form, $\omega_k^L$.
\label{lem:stormer}\end{lem}
\begin{pf}
Consider a discrete trajectory $(q_k)_k$ that verifies equation (\ref{eq:stormer_l}). Then we have:
\begin{eqnarray}
d S_d^L&=&\tau\sum_{k=1}^{n-1} \left(\pfrac{L_d(\qk,\Dt \qk)}{\qk^i}-\Dt \pfrac{L_d(\qkmo,\Dt
\qkmo)}{\Dtd \qkmo^i }\right)d\qk^i\nonumber\\&& +\Dt \left(\pfrac{ L_d(\qkmo,\Dt \qkmo)}{\Dt \qk^i}d
\qk^i\right)\,.
\label{eq:tmp1}
\end{eqnarray}
Since the $q_k$'s verify equation (\ref{eq:stormer_l}), and $d^2=0$, equation (\ref{eq:tmp1})
yields:
\begin{equation}
\begin{array}{ccc}
d(\Dt\theta_k^L)=0\,, & \textrm{that is,}&\omega_{k+1}^L=\omega_k^L\,.
\end{array}
\label{eq:stormer_tmp}
\end{equation}
We conclude that $\omega^L_k$ is preserved along the discrete trajectory
\end{pf}
As we mentioned earlier, because DVPII acts on the tangent bundle it provides results very similar
to the continuous case as attested by the form of the Lagrangian $2$-form. This is to be compared
with the Lagrangian two-form arising in the continuous case:
 \begin{equation}
\omega^L=\ppfrac{L}{q^i}{\dot{q}^j}d q^i\wedge dq^j +\ppfrac{L}{\dot{q}^i}{\dot{q}^j}d \dot{q}^i\wedge
dq^j\,.
 \end{equation}

Note that conservation of the Lagrangian two-form is a consequence of using the Leibnitz law, and
therefore does not depend on the definition of the discrete Lagrangian. In the remainder of this
section we use different discrete Lagrangian functions, but the same Leibnitz law. Thus 
lemma \ref{lem:stormer} still applies.

More generally, we can derive St\"{o}rmer's rule using $$L_d(\qk,\Dt\qk)=\lambda
L(\qk,\Dt\qk)+(1-\lambda)L(\qk+\tau\Dt\qk,\Dt\qk)\,,$$ for any $\lambda$ in $\mathbb{R}$.  A
particular case of interest is $\lambda=\undemi$ which yields a symmetric version of St\"{o}rmer's rule
also called the velocity Verlet method\cite{swo82}. For this value of $\lambda$, we define the
associated discrete momenta using the Legendre transform (equation (\ref{eq:legendre})):
\begin{eqnarray}
\pko&=&\pdk\\
    &=&D_2 L_d(\qk,\Dt\qk)\\
    &=&M \Dt \qk-\undemi \tau\nabla V(\qk+\Dt \qk)\,,
\label{eq:legendre_vv}
\end{eqnarray}
that is:
\begin{equation}
\qko=\qk+ \tau M^{-1}(\pko+\undemi \tau\nabla V(\qko))\,.
\label{eq:l_vv1}
\end{equation}
Moreover, from equation (\ref{eq:stormer_del}) we obtain:
\begin{equation}
\pko=\pk+\tau\frac{-\nabla V(\qk)-\nabla V(\qko)}{2}\,.
\label{eq:l_vv2}
\end{equation}
Equations (\ref{eq:l_vv1}) and (\ref{eq:l_vv2})
 define the velocity Verlet algorithm.

We now focus on the Newmark algorithm which is usually written for the system $L=\undemi
\dot{q}^T M \dot{q}-V(q)$ as a map given by $(\qk,\dot{q}_k)\mapsto (\qko,\dot{q}_{k+1})$ satisfying the
 implicit relations:
\begin{eqnarray}
\qko&=&\qk+ \tau \dot{q}_k+\frac{\tau^2}{2}[(1-2\beta)a_k+2 \beta a_{k+1}]\label{eq:newmark1}\,,\\
\dot{q}_{k+1}&=&\dot{q}_k+\tau[(1-\gamma)a_k+\gamma a_{k+1}]\label{eq:newmark2}\,,\\
a_k&=&M^{-1}(-\nabla V(\qk))\,,
\end{eqnarray}
where the parameters $\gamma \in [0,1]$ and $\beta \in [0,\undemi]$. For $\gamma=\undemi$ and any
$\beta$  the Newmark algorithm can be generated from DVPII as a particular case of the St{\"o}rmer
rule where $\qdk$ and $L_d$ are chosen as follows:
$$\qdk=\qk-\beta \tau^2 a_k\,,$$
and
$$L_d(\qdk,\Dtd \qdk)=\undemi \dot{\qdk}^T M
\dot{\qdk}-\tilde{V}(\qdk)\,,$$
with $\tilde{V}$, the modified potential, satisfying $\nabla\tilde{V}(\qdk)=\nabla V(\qk)$. Since the
derivative operator is the same as above, the discrete Hamilton's principle yields St{\"o}rmer's
equation where $\qk$ is replaced by $\qdk$, that is:
\begin{equation}
q^d_{k+1}=2\qdk-q^d_{k-1}+\tau^2 M^{-1}(-\nabla \tilde{V}(\qdk))\,.
\label{eq:stormer_newmark}
\end{equation}
Equation (\ref{eq:stormer_newmark}) simplifies to
\begin{equation}
\qko-2\qk+\qkmo=\tau^2(\beta a_{k+2}+(1-2 \beta)a_{k+1}+\beta a_{k-1})\,.
\end{equation}
This last equation corresponds to the Newmark algorithm for the case $\gamma=\undemi$. Lemma
\ref{lem:stormer} guarantees that  the Lagrangian two-form $$\omega_k^L=d(D_2 L_d(\qdk,\Dtd
\qdk)dq^d_{k+1})$$ is preserved along the discrete trajectory.

\subsubsection{From the Hamiltonian point of view}
The St{\"o}rmer, velocity Verlet, and Newmark algorithms can also be derived using a phase space
approach, i.e., the DMHP. For St\"{o}rmer's rule, the Legendre transform yields:
\begin{equation}
\pko=M \Dt \qk\,.
\label{eq:legendre_stormer}
\end{equation}
The discrete Hamiltonian function is defined from equation
(\ref{eq:discrete_ham}):
\begin{equation}
H_d(\qk,\pko)=\undemi \pko^T M^{-1}\pko +V(\qk)\,,
\end{equation}
and discrete equations of motion are obtained from the DMHP\footnote{$\qdk=\qk$ and $\pdk=\pko$}
(theorem (\ref{def:discrete_pcp_la})). We skip a few steps in the evaluation of the variations of
$S^H_d$ to finally find:
\begin{eqnarray}
\delta S_d^H&=&\delta \left(\tau\sum_{k=0}^{n-1} \langle\pko, \Dt\qk\rangle -H_d(\qk,\pko)\right)\\
&=&\tau\sum_{k=0}^{n-1} \langle\Dt \qk -D_2H_d(\qk,\pko),\delta \pko\rangle-\langle\Dt \pk
+D_1H_d(\qk,\pko),
\delta \qk\rangle \nonumber\\
&&\qquad+\langle p_{n}, \delta \qn\rangle-\langle p_0, \delta
\qz\rangle\,.
\end{eqnarray}
If we impose the variations of the action $S^H_d$ to be zero for any $(\delta \qk,\delta \pko)$ and
$\delta q_0=\delta q_n=0$, we obtain:
\begin{eqnarray}
\Dt \qk&=&\pko\label{eq:stormer_h1}\,,\\
\Dt\pk&=&-\nabla V(\qk)\,.\label{eq:stormer_h2}
\end{eqnarray}
Elimination of the $\pk$'s yields St{\"o}rmer's rule.

To recover the velocity Verlet scheme from the Hamiltonian point of view, one needs to solve for
$\Dt\qk$  as a function of $(\qk,\pko)$ in equation (\ref{eq:legendre_vv}). Suppose this has been
done and that $\Dt\qk=f(\qk,\pko)$, then
\begin{equation}
H_d(\qk,\pko)=\langle \pko,f(\qk,\pko)\rangle - L_d(\qk,f(\qk,\pko))
\label{eq:ham_vv}\,.
\end{equation}
Taking the variation of the action $S^H_d$ yields the following discrete Hamilton's equations:
\begin{eqnarray}
\Dt\qk&=&D_2 H_d(\qk,\pko)\label{eq:ham_vv1}\,,\\
\Dt\pk&=&-D_1 H_d(\qk,\pko)\,.
\label{eq:ham_vv2}
\end{eqnarray}
On the other hand, equation (\ref{eq:ham_vv}) provides the following relationships:
\begin{multline}
D_1 H_d(\qk,\pko)=D_1 f(\qk,\pko) ( \pko- D_2 L_d(\qk,f(\qk,\pko)))\\
-D_1 L_d(\qk,f(\qk,\pko))\,,\label{eq:ham_vv3}
\end{multline}
\begin{equation}
D_2 H_d(\qk,\pko)=\Dtd\qk +D_2 f(\qk,\pko) (\pko-D_2 L_d(\qk,f(\qk,\pko)))\,.
\label{eq:ham_vv4}
\end{equation}
Combining equations (\ref{eq:ham_vv1}) and (\ref{eq:ham_vv2}) together with equations
(\ref{eq:ham_vv3}) and (\ref{eq:ham_vv4}) yields the Velocity Verlet algorithm (equations
(\ref{eq:l_vv1}) and (\ref{eq:l_vv2})).

We now prove that the scheme we obtained is symplectic. As in the Lagrangian case, the proof
differs from the usual one that consists in computing $d\pko \wedge d\qko$, in that it relies on
fundamental properties of DVPII and on the use of the Leibnitz law.
\begin{lem}
The algorithm defined by equations (\ref{eq:ham_vv1})-(\ref{eq:ham_vv2}) is symplectic.
\end{lem}
\begin{pf}
We have:
\begin{eqnarray}
d S_d^H&=&d \left(\tau\sum_{k=0}^{n-1} \langle\pko, \Dt\qk\rangle -H_d(\qk,\pko)\right)\,,\\
&=&\tau\sum_{k=0}^{n-1} \langle\Dt \qk -D_2H_d(\qk,\pko),d \pko\rangle-\langle\Dt \pk
+D_1H_d(\qk,\pko),d\qk\rangle \nonumber\\
&&\,{+}\Dt \langle p_k ,dq_k\rangle\,.
\end{eqnarray}
Hence, since $(q_k,p_k)$ verifies equations (\ref{eq:ham_vv1})-(\ref{eq:ham_vv2}) and $d^2=0$, we
obtain:
\begin{equation}
\Dt (d\pk\wedge d\qk)=0\,.
\end{equation}
The symplectic two-form $d\pk\wedge d\qk$ is preserved along the trajectory.
\end{pf}

Finally, we can also derive the Newmark methods from the Hamiltonian point of view. The Legendre
transform yields:
\begin{equation}
\pdk=\pfrac{L_d(\qdk,\Dtd \qdk)}{\Dtd \qdk}=M \Dtd \qdk\,.
\end{equation}
The Newmark algorithm is again a particular case of the St{\"o}rmer rule where $(\qk,\pko)$ is replaced
by $(\qdk,\pdk)$:
\begin{eqnarray}
\Dtd \qdk&=&\pdk\,,\\
\Dtd \pdk&=&-\nabla \tilde{V}(\qdk)\,.
\end{eqnarray}
Defining $\dot{q}_k$ from $\pk$ as $$\dot{q}_k=M^{-1}\pdk+\frac{\tau}{2}a_k$$ allows us to recover
the Newmark scheme for $\gamma=\undemi$ (equations (\ref{eq:newmark1}) and (\ref{eq:newmark2})).
From the above lemma, we obtain that the symplectic two-form $d\pdk
\wedge d\qdk$ is preserved along the trajectory.

\subsection{Midpoint rule}
The midpoint rule has been extensively studied and a complete study of its properties can be found
in the literature. It is a particular case of the Runge-Kutta algorithm, but can also be derived as
a variational integrator (see for instance
\cite{wu89,shi92,mar01}). The derivation of this scheme has been done by Guo, Li and Wu \cite{guo02li3} for the
Hamiltonian point of view. In the next section we present the Lagrangian point of view and then
recall the Guo, Li and Wu main results, the goal of this section being to illustrate the use of DVPII
with other discretization and discrete derivative operator.

\subsubsection{From the Lagrangian point of view}

Given a Lagrangian $L(q,\dot{q})$, define the discrete Lagrangian by:
\begin{equation}
L_d(\qdk,\Dtd\qdk)=L(\qdk,\Dtd\qdk)\,,
\end{equation}
where $\qdk=\frac{\qko+\qk}{2}$, and $\Dtd=R_{\tau/2}-R_{-\tau/2}$ where the operator $R_\tau$ is
the translation by $\tau$. One can readily verify that
 $\Dtd \qdk=\Dt \qk$ and that $\Dtd$ verifies
the usual Leibnitz law:
\begin{eqnarray}
\Dtd (\fdk \gdk)=\Dtd \fdk \cdot \gdk+\fdk \cdot\Dtd \gdk\,,
\label{eq:leibnitz_midpoint}
\end{eqnarray}
where $f_k=f(t_k)$ and $g_k=g(t_k)$ are functions of time and $\fdk=\frac{f_{k+1}+f_k}{2}$.
Applying the discrete Hamilton's principle yields:
\begin{eqnarray}
\delta S^L_d&=&\tau\sum_{k=0}^{n-1} \delta L_d(\qdk,\Dtd \qdk)\\
&=&\tau\sum_{k=0}^{n-1} \langle D_1 L_d(\qdk,\Dtd \qdk),\delta \qdk\rangle +\langle D_2
L_d(\qdk,\Dtd \qdk),\delta \Dtd \qdk\rangle \,.\label{eq:midpoint_temp0}
\end{eqnarray}
From the Legendre transform (equation (\ref{eq:legendre})), we define the associated momentum:
\begin{equation}
\frac{\pko+\pk}{2}=\pdk=D_2 L_d(\qdk,\Dtd \qdk)\,.
\end{equation}
Then, equation (\ref{eq:midpoint_temp0}) becomes:
\begin{eqnarray}
\delta S^L_d&=&\tau\sum_{k=0}^{n-1} \langle D_1 L_d(\qdk,\Dtd \qdk),\delta \qdk\rangle +\langle \pdk,
 \delta \Dtd \qdk\rangle\\
&=&\tau\sum_{k=0}^{n-1} \langle D_1 L_d(\qdk,\Dtd \qdk)-\Dtd \pdk,\delta \qdk\rangle+ \langle
\pn ,\delta \qn\rangle -\langle \pz,\delta \qz \rangle\,.
\end{eqnarray}
If we require the variations of the action to be zero for any choice of $\delta \qdk,\:
k\in[1,n-1]$, and $\delta
\qz=\delta \qn=0$, we obtain discrete Euler-Lagrange equations for the midpoint scheme:
\begin{eqnarray}
\frac{\pko-\pk}{h}&=&\Dtd \pdk\nonumber\\
&=&D_1 L_d(\qdk,\Dtd \qdk) \nonumber\\
&=&D_1 L_d(\frac{\qko+\qk}{2},\frac{\qko-\qk}{h})
\label{eq:midpoint_l1}\,,\\
\frac{\pko+\pk}{2}&=&\pdk\nonumber\\
&=&D_2 L_d(\qdk,\Dtd \qdk)\nonumber\\
&=&D_2 L_d(\frac{\qko+\qk}{2},\frac{\qko-\qk}{h})\,.
\label{eq:midpoint_l2}
\end{eqnarray}

\begin{lem}
The midpoint scheme defines a symplectic algorithm.
\end{lem}
\begin{pf}
The proof proceeds as for the St{\"o}rmer rule:
\begin{eqnarray}
dS^L_d&=&\tau\sum_{k=0}^{n-1} \langle D_1 L_d(\qdk,\Dtd \qdk),d \qdk\rangle +
\langle \pdk, d \Dtd \qdk\rangle\\
&=&\tau\sum_{k=0}^{n-1} \langle D_1 L_d(\qdk,\Dtd \qdk)-\Dtd \pdk,d \qdk\rangle+ \Dtd \langle\pdk,
d\qdk\rangle\,.
\end{eqnarray}
Since $d^2=0$ and $(\qk,\pk)$ verifies equations (\ref{eq:midpoint_l1})-(\ref{eq:midpoint_l2}), we
obtain:
\begin{equation}
\Dtd (d\pdk \wedge d\qdk)=0\,.
\end{equation}
A straight forward computation shows that $\Dtd (d\pdk \wedge d\qdk)=\Dt (d\pk \wedge d\qk)$, i.e.,
the symplectic two-form $\omega_k=d\pk \wedge d\qk$ is preserved along the trajectory.
\end{pf}

\subsubsection{From the Hamiltonian point of view}
Let $H_d(\qdk,\pdk)=H(\qdk,\pdk)$ or equivalently define $H_d$ from $L_d$ via equation
(\ref{eq:discrete_ham}) and let $(\qdk,\pdk)=(\frac{\qko+\qk}{2},\frac{\pko+\pk}{2})$. Then the
DMHP (\ref{def:discrete_pcp_la}) yields:
\begin{eqnarray}
\frac{\qko-\qk}{h}&=&\Dtd\pdk\nonumber\\
&=&D_2H_d(\qdk,\pdk)\nonumber\\
&=&\pfrac{H}{p}(\frac{\qko+\qk}{2},\frac{\pko+\pk}{2})\,,\\
\frac{\pko-\pk}{h}&=&\Dtd\pdk\nonumber\\
&=&-D_1H_d(\qdk,\pdk)\nonumber\\
&=&-\pfrac{H}{q}(\frac{\qko+\qk}{2},\frac{\pko+\pk}{2})\,.
\end{eqnarray}

\begin{lem}
The midpoint scheme defines a symplectic algorithm.
\end{lem}
\begin{pf}
The proof is straightforward. We compute $d^2 S_d^H$ assuming $(\qk,\pk)$ verifies the above
equations of motion.
\end{pf}

To conclude, we have illustrated the use of the discrete variational principles
(definitions \eqref{defi:discrete_ham_pcp} and \eqref{def:discrete_pcp_la}) and derived discrete equations of
motion. One can readily verify that both variational principles yield the same discrete equations,
as in the continuous case.
%
Other schemes can be recovered in the same way, and we do not know yet if all classical symplectic
algorithms can be derived from DVPII.
For instance, we have been able to recover the conditions for the partitioned Runge-Kutta algorithm
to be symplectic from the Lagrangian point of view but so far it is not clear to us how it can be
done using the Hamiltonian approach (definition \eqref{def:discrete_pcp_la}).

\subsection{Numerical example}
Symplectic integrators are usually used as numerical integrators that preserve the qualitative
behavior of dynamical systems and are especially valuable for long time simulations. However, these
are not the only uses of symplectic integrators. In this section we present an
aspect of symplectic integrators that we have not seen pointed out in the literature: we show that
they allow one to recover the generating functions for the phase flow canonical transformation,
whereas numerical integrators do not, even over a short period of time (applications of this result
can be found in \cite{gui04aiaa}).

Let us first recall two results from the Hamilton-Jacobi theory.
\begin{prop}
The transformation induced by the phase flow is canonical.
\end{prop}
\begin{prop}
Let $(P_1,\omega_1)$ and $(P_2,\omega_2)$ be symplectic manifolds, $\pi_i:P_1\times P_2
\rightarrow P_i$ the projection onto $P_i$, $i=1,2$, and
\begin{equation}
\Omega=\pi_1^*\omega_1-\pi_2^*\omega_2\,.
\end{equation}
Then:
\begin{enumerate}
\item $\Omega$ is a symplectic form on $P_1\times P_2$;
\item a map $f: P_1 \rightarrow P_2$ is symplectic if and only if
$i_f^*\Omega=0$, where $i_f: \Gamma_f\rightarrow P_1\times P_2$ is inclusion and $\Gamma_f$ is the
graph of $f$.
\end{enumerate}
\end{prop}
Hence, by the Poincar{\'e} lemma, if $f$ is canonical there exists a function $S$ such that
\begin{equation}
i_f^*\Theta=dS\,,
\label{eq:can_gen}
\end{equation}
where $\Omega=-d\Theta$. $S$ is called a generating function. If $(q^i,p_i)$ are coordinates on
$P_1$ and $(Q^i, P_i)$ are coordinates on $P_2$, then $\Gamma_f$ can be endowed with a chart in
several ways. For instance, $S$ may appear as a function of $(q^i, Q^i)$ or of $(q^i, P_i)$, and so
forth depending of the choice of $\Theta$. Let $\theta_1=p_idq^i$ and $\theta_2=P_idQ^i$, then
$i_f^*\Theta=i_f^*\pi_1^*\theta_1-i_f^*\pi_2^*\theta_2=(\pi_1 \circ i_f)^*p_idq^i-(\pi_2 \circ
i_f)^*P_idQ^i$. In this case, $S$ is a function of $(q^1,\cdots,q^n,, Q^1,\cdots,Q^n)$. From
$$dS= \pfrac{S}{q^i}dq^i+\pfrac{S}{Q^i}dQ^i\,,$$ we conclude,
using equation (\ref{eq:can_gen}) that:
\begin{equation} p_i=\pfrac{S}{q_i} \qquad
P_i=-\pfrac{S}{Q^i}
\label{eq:can_gen2}
\end{equation}

Suppose that $f$ is the phase flow $\Phi$, then equation (\ref{eq:can_gen2}) defines a relationship
between flow and the gradient of the generating function. In particular, if the generating function
$S(q,q_0,t)$ exists and the flow is defined as follows:
\begin{equation}
\Phi : (q_0,p_0,t)\mapsto (q(t),p(t))=(\Phi_1(q_0,p_0),\Phi_2(q_0,p_0))\,,
\end{equation}
then, from the local inverse function theorem\footnote{$|\pfrac{\Phi}{p_0}|\ne 0$ since we assume that $S$
exists}, there exist two functions $S_1$ and $S_2$ such that:
\begin{eqnarray}
p_0&=&S_1(q,q_0,t)\,,\label{eq:s1}\\
p&=&\Phi_2(q_0,S_1(q,q_0,t))\equiv S_2(q,q_0,t)\,.\label{eq:s2}
\end{eqnarray}
From equation (\ref{eq:can_gen2}), we conclude that $S_1$ and $S_2$ are the gradient of $S$ and
therefore should verify\footnote{Since their exists an open set on which the generating functions
are smooth, Schwartz's theorem yields $\ppfrac{S}{q_0}{q}=\ppfrac{S}{q}{q_0}$.}:
\begin{equation}
\ppfrac{S}{q_0}{q}\equiv\pfrac{S_1}{q}(q,q_0,t)=\pfrac{S_2}{q_0}(q,q_0,t)\equiv \ppfrac{S}{q}{q_0}\,.
\label{eq:exactness}
\end{equation}
Since symplectic integrators preserve the symplectic two-form, the exactness condition (equation
(\ref{eq:exactness})) is satisfied whereas it is not using numerical integrators.

\subsubsection{Harmonic Oscillator}
We start with a trivial example, the harmonic oscillator, because its study allows us to introduce
techniques and discuss issues that arise in the next more sophisticated example. The Hamiltonian
function for the harmonic oscillator is quadratic:
\begin{equation}
H(q,p)=\frac{1}{2m} p^2+\frac{k}{2} q^2\,.
\end{equation}
It is a linear system so the phase flow is also linear:
\begin{eqnarray}
\Phi_1(q_0,p_0)&=&a_{11}(t) q_0+a_{12}(t)p_0\,\\
\Phi_2(q_0,p_0)&=&a_{21}(t) q_0+a_{22}(t)p_0\,.
\end{eqnarray}
Substituting these expressions into Hamilton's equations and balancing terms of the same order
yield:
\begin{equation}
\left\{\begin{array}{rcl}
\dot{a}_{11}(t)&=&a_{21}(t)/m\\
\dot{a}_{12}(t)&=&a_{22}(t)/m\\
\dot{a}_{21}(t)&=&k a_{11}(t)\\
\dot{a}_{22}(t)&=&k a_{12}(t)
\end{array}\right.
\label{eq:stm_harm}
\end{equation}

\begin{figure*}[htb]
\begin{center}
\subfigure[Midpoint scheme with fixed time step $\tau=0.01$]{\includegraphics[width=0.45\linewidth]{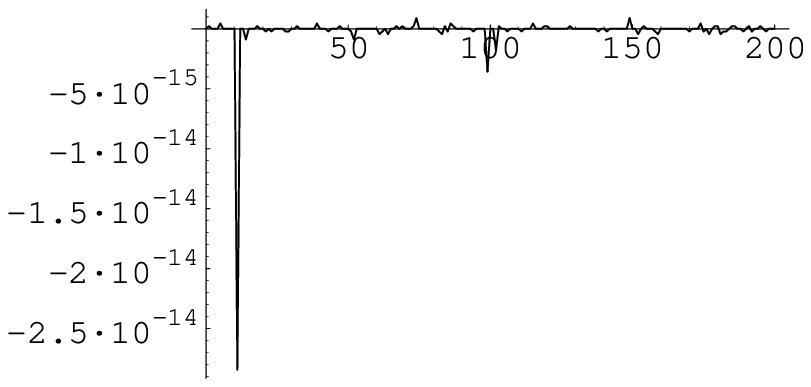}}
\subfigure[Implicit Gauss Runge-Kutta algorithm of order $8$]{\includegraphics[width=0.45\linewidth]{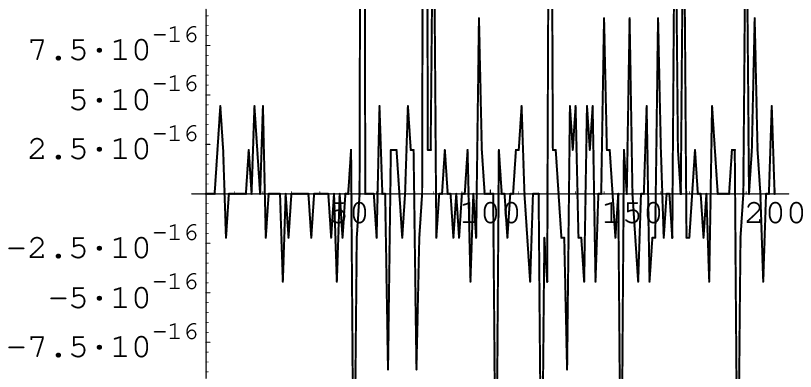}}
\subfigure[Explicit Runge-Kutta algorithm of order $8$]{\includegraphics[width=0.45\linewidth]{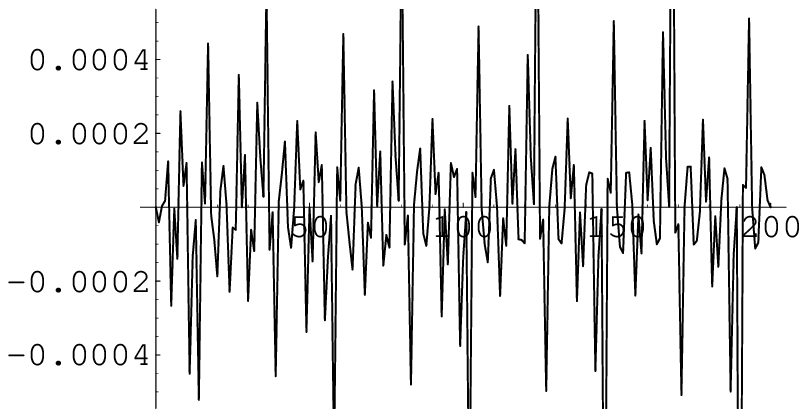}}
\caption{\label{fig:harmonic}Exactness condition using $3$ different integrators}
\end{center}
\end{figure*}

In figure \ref{fig:harmonic}, we plot $\Delta=\pfrac{S_1}{q}(q,q_0,t)-\pfrac{S_2}{q_0}(q,q_0,t)$
over the time interval $[0,100]$ using the midpoint scheme with fixed time step, a symplectic Gauss
implicit Runge-Kutta algorithm of order $8$ with fixed time step and a non symplectic Runge-Kutta
integrator of order $8$ to integrate equations (\ref{eq:stm_harm}). We remark that only
symplectic integrators allow us to recover the generating functions because the exactness condition
is exactly verified. We point out that even over a short time span, numerical integrators fail to
satisfy the exactness condition.

\subsubsection{Earth orbit}
This example was first encountered  by V.M. Guibout and D.J. Scheeres \cite{gui03,gui04aiaa}
while studying spacecraft formation flight. Consider an orbital problem about the Earth modelled by
a non-spherical body (we take into account $J_2$ and $J_3$ gravity coefficients). The Hamiltonian
of the system is given by
\begin{eqnarray}
H &=& \undemi (p_x^2+p_y^2+p_z^2) \nonumber\\
&&- \frac{1}{\sqrt{x^2+y^2+z^2}} \left( 1 -\frac{R^2}{2 r_0^2 (x^2+y^2+z^2)}
(3\frac{ z^2}{x^2+y^2+z^2}-1)J_2\right.\nonumber\\
&& \left. -\frac{R^3}{2 r_0^3 (x^2+y^2+z^2)^2}(5\frac{ z^3}{x^2+y^2+z^2}-3 z)J_3
\right)\,,
\end{eqnarray}
where $$\begin{array}{cc} GM=398,600.4405\; km^3s^{-2}\,,&
R =6,378.137\; km\,,\\
J2=1.082626675\cdot 10^{-3}\,,& J3=2.532436\cdot 10^{-6}\,,
\end{array}$$
and all the variables are normalized ($r_0$ is the initial radius of a trajectory):
\begin{equation}
\begin{disarray}{rclccccccccccc}
&&&x&\rightarrow&x r_0 \,,&\quad & y&\rightarrow&y r_0\,,&\quad& z&\rightarrow&z r_0\,,\\
t&\rightarrow&t \sqrt{\frac{r_0^3}{GM}}\,,& p_x&\rightarrow&p_x
\sqrt{\frac{GM}{r_0}}\,,&\quad&p_y&\rightarrow&p_y
\sqrt{\frac{GM}{r_0}}\,,&\quad& p_z&\rightarrow&p_z  \sqrt{\frac{GM}{r_0}}\,.
\end{disarray}
\end{equation}
We choose the nominal trajectory to be a highly eccentric orbit. The initial conditions for the
nominal trajectory are in normalized units ( $r_0=7000 \;km$):
\begin{equation}
\begin{array}{ccccccccccc}
x &=& 1\,,&\quad& y &=& 0\,,&\quad&
z &=& 0\,,\\
p_x&=& 0 \,,&\quad& p_y&=&  \sqrt{\frac{13}{10}}\cos(\pi/3)\,,&\quad& p_z&=&
v\sqrt{\frac{13}{10}}\sin(\pi/3)\,.
\end{array}
\end{equation}
At the initial time, $e=0.3$, $i=\frac{\pi}{3} \;rad$, $\omega=0$ and $\Omega=0$.

 This system is non-integrable and has non-trivial dynamics. The phase
flow is not known globally but techniques have been developed to evaluate it locally (see Guibout
and Scheeres
\cite{gui04}).

Consider a given trajectory called the nominal trajectory $(q^0(t),p^0(t))$, then the dynamics of
the relative motion of a particle about this trajectory is Hamiltonian and is described by the
following Hamiltonian function $H^h(X^h,t)$:
\begin{equation}
\sum_{p=2}^\infty \sum_{\substack{i_1,\cdots,i_{2n}=0\\
\sum_{k=1}^{2n} i_k=p}
}^p
\frac{1}{i_1!\cdots i_{2n}!}
\frac{\partial^{p} H}{\partial q_1^{i_1}\cdots\partial q_n^{i_n}
\partial p_1^{i_{n+1}}\cdots\partial p_n^{i_{2n}}}(q^0,p^0,t)
{X^h_1}^{i_1}\dots{X^h_{2n}}^{i_{2n}}\,,
\label{eq:l_hamilton_n}
\end{equation}
where $X^h$ is the relative state vector $(\hq,\hp)$.

In the same way, we expand in Taylor series the phase flow for the relative motion, and substitute
its expression into Hamilton's equations for $H^h$. When studying spacecraft formation flight we
often assume that the spacecraft stay close to each other and therefore, one may approximate the
dynamics of the formation by truncating the above Taylor series. Suppose we keep only terms of
order less that $N$. Then balancing terms of the same order in Hamilton's equations yields a set of
ordinary differential equations (the procedure is the same as in the harmonic oscillator example
but here there are non linear terms up to order $N$). We use the midpoint scheme with fixed-time
step ($\tau=0.01$), a symplectic Gauss implicit Runge-Kutta algorithm of order $4$ with fixed time
step ($\tau=0.01$) and $Mathematica^\copyright$ built-in numerical integrator
\textit{NDSolve}\footnote{\textit{NDSolve} switches between a non-stiff Adams method and a stiff
Gear method.} to integrate the flow up to order $N=4$. Once the Taylor series of the flow is known,
we find $S_1$ and $S_2$ by a series inversion.
Then we check the exactness conditions defined by equation (\ref{eq:exactness}) (there are several terms involved since we are dealing with a nonlinear system of dimension $6$ ). We find that after $10 \pi$ units of time, $\|\pfrac{S_1}{q}-\pfrac{S_2}{q_0}\|\le \eta$, where $\eta=10^{-11}$ using the midpoint scheme, $\eta=10^{-11}$ with the symplectic implicit Gauss Runge-Kutta algorithm and $\eta=10^{-3}$ with the built-in function \textit{NDSolve}. Again, only symplectic algorithms allow us to  recover the generating functions.

\section{Energy conservation}
\label{sec:energy}
Symplectic integrators do not conserve energy and in general induce bounded energy error. There are several works on analyzing the energy error, we refer to Hairer and Lubich \cite{hai01lub} and Hairer, Lubich and Wanner \cite{hai02lub} and references therein for more details. In this section, we enhance DVPII so that energy conservation is imposed. By considering the time as a coordinate and by adding an independent parameter $\tau$, DVPII yields symplectic energy conserving algorithms. For certain problems, such algorithms may provide better performance\footnote{To quantify the performance of an algorithm, not only we look at its accuracy but we also evaluate its ability to predict the qualitative behavior of the system. In that sense, symplectic-energy conserving algorithms may not predict qualitative behavior better that symplectic algorithms.}, but the contrary may also happen \cite{hai00lub,sim93gon}. The method we develop in this section is variational and allows us to recover Shibberu's algorithm \cite{shi92} for Hamiltonian systems and is equivalent to the Kane, Marsden and Ortiz \cite{kan99} method for Lagrangian systems.

\subsection{Generalized variational principles}

\subsubsection{Generalized Hamilton's principle}
Let us first recall Hamilton's principle for dynamical systems for which time is considered as a generalized coordinate. Such a formulation is typically used in relativity where the time coordinate is equivalent to the space coordinates.

Consider a Lagrangian $L(q,\dot{q})$ and define the \textit{parametric} Lagrangian $$\bar{L}(q,t,q',t')= t' L(q,\frac{q'}{t'},t)\,,$$ where $'=\frac{d}{d\tau}$ and $\tau$ is an independent parameter that parameterizes the trajectory and the time. Then the generalized Hamilton's principle reads:
\begin{defn}
Critical points of $\int_{t_0}^{t_f} \bar{L}(q,\frac{q'}{t'},t)d\tau$ in the class of curves
$(q(\tau),t(\tau))$ with endpoints $(q_0,t_0)$ and $(q_f,t_f)$ correspond to trajectories of the Lagrangian systems going from $(q_0,t_0)$ to $(q_f,t_f)$.
\end{defn}
The generalized Hamilton's principle yields the following set of equations:
\begin{eqnarray}
\pfrac{\bar{L}}{t}-\frac{d}{d\tau}\pfrac{\bar{L}}{t'}&=&0\,,\\
\pfrac{\bar{L}}{q}-\frac{d}{d\tau}\pfrac{\bar{L}}{q'}&=&0\,.
\end{eqnarray}
Replacing the parametric Lagrangian by the Lagrangian of the system simplifies the above equations
to:
\begin{eqnarray}
t'\pfrac{L}{t}-\frac{d}{d\tau}L+\frac{d}{d\tau}\left(\pfrac{L}{\dot{q}}\frac{q'}{t'}\right)&=&0\,,\label{eq:gel_eq1}\\
t'\pfrac{L}{q}-\frac{d}{d\tau}\pfrac{L}{\dot{q}}&=&0\,.\label{eq:gel_eq2}
\end{eqnarray}
These $n+1$ equations should be compared to the $n$ equations obtained when the trajectory is
parameterized by the time:
\begin{equation}
\pfrac{L}{q}-\frac{d}{dt}\pfrac{L}{\dot{q}}=0\,.\label{eq:el_eq}
\end{equation}
Since $\frac{d}{d\tau}=t'\frac{d}{dt}$, we conclude that the space components of the generalized
Euler-Lagrange equations (equation (\ref{eq:gel_eq2})) are a multiple by $t'$ of the original
Euler-Lagrange equations (equation (\ref{eq:el_eq})). Also, their
 time component (equation
(\ref{eq:gel_eq1})) is a linear combination of the
components of equation (\ref{eq:el_eq}) (the
 sum of each component multiplied by
$q'$). All $n+1$ generalized Euler-Lagrange equations are thus consistent with the original
equations but there is no unique solution because they are satisfied by any parameterization. To
get a unique solution, it is necessary to add to the generalized Hamilton's principle an additional
condition fixing the parameterization. As we will see in the next section, in discrete settings we
do not have this freedom anymore. The discrete counter-part of equation (\ref{eq:gel_eq1})
corresponds to an energy constraint that fully specifies the time parameterization, i.e., the time
step.

\subsubsection{Generalized discrete Hamilton's principle (GDHM)}
In contrast with the variational principles introduced in the first part of this paper, we do not
set the time step,  i.e., we let the time act as a variable by adding an independent parameter
$\tau_k$ such that $t_k=t(\tau_k)$ and $\tau_{k+1}-\tau_k=\tau$, $\tau$ being a constant.
$t_k$ is now a coordinate that plays the same role as $q_k$, $M_k$ is the extended configuration
space $(\qk,t_k)$, ${\mathcal M}=\bigcup M_k$ and ${\mathcal T}=\{(\tau_k)_{k\in [1,n]}\}$. Define
the modified discrete Lagrangian $\bar{L}_d$:
\begin{equation}
\bar{L}_d(\qdk,\tdk,\Dtd \qdk,\Dtd \tdk)=\Dtd \tdk L_d(\qdk,\frac{\Dtd \qdk}{\Dtd \tdk},\tdk)\,,
\end{equation}
where $L_d$ is the discrete Lagrangian previously defined. In addition, since we are interested in
conservation of energy we only consider system that are time independent. As a consequence, $L_d$
does not depend on time and $\pfrac{\bar{L}_d}{\tdk}=0$.

\begin{defn}[Generalized Discrete Hamilton's Principle (GDHP)]
\label{defi:m_discrete_ham_pcp}
Critical points of the discrete action
\begin{equation}
S^L_d=\sum_{k=0}^{n-1} \bar{L}_d(\qdk,\tdk,\Dtd \qdk,\Dtd \tdk) \tau\,,
\end{equation}
in the class of discrete curves $(\qdk,\tdk)_k$ with endpoints $(\tau_0,t_0,q_0)$ and
$(\tau_n,t_n,q_n)$ correspond to trajectories of the discrete Hamiltonian system going from $(t_0,q_0)$ to $(t_n,q_n)$:
\end{defn}

Again, to proceed to the derivation of the equations of motion we need to specify the derivative
operator.

\subsubsection{Generalized discrete modified Hamilton's principle}

\begin{defn}
Let $\bar{L}_d$ be a discrete Lagrangian on $T{\mathcal M}$ and define the discrete Legendre
transform (or discrete fiber derivative) $\mathbb{F}L:T{\mathcal M}\rightarrow T^*{\mathcal M}$
which maps the discrete extended phase space $T{\mathcal M}$ to $T^*{\mathcal M}$ by
\begin{equation}
(\qdk,\tdk,\Dtd \qk,\Dtd \tdk)\mapsto(\qdk,\tdk,\pdk,\edk)\,,
\end{equation}
where
\begin{equation}
\pdk=\pfrac{\bar{L}_d(\qdk,\tdk,\Dtd \qdk,\Dtd \tdk)}{\Dtd \qdk}\,,\quad
\edk=\pfrac{\bar{L}_d(\qdk,\tdk,\Dtd \qdk,\Dtd \tdk)}{\Dtd \tdk}\,.
\label{eq:m_legendre}
\end{equation}
The Legendre transform as defined by equations (\ref{eq:m_legendre}) is equivalent to the previous
definition (equation \eqref{eq:legendre}). Indeed,
$$\pfrac{\bar{L}_d(\qdk,\tdk,\Dtd \qdk,\Dtd \tdk)}{\Dtd \qdk}=\pfrac{L_d(\qdk,\frac{\Dtd
\qdk}{\Dtd \tdk})}{\Dtd \qdk}=D_2L_d(\qdk,\Delta_t^d \qdk)\,,$$
where $\Delta_t^d = \frac{\Dtd}{\Dtd \tdk}$ represent the discrete derivative with respect to time.
\end{defn}
If the discrete fiber derivative is a local isomorphism, $\bar{L}_d$ is called regular and if it is
a global isomorphism we say that $\bar{L}_d$ is hyperregular. If $\bar{L}_d$ is hyperregular, we
define the corresponding discrete Hamiltonian function on $T^*{\mathcal M}$ by
\begin{equation}
\bar{H}_d(\qdk,\tdk,\pdk,\edk)=\langle\pdk, \Dtd \qdk\rangle - \bar{L}_d(\qdk,\Dtd\qdk)\,,
\label{eq:m_discrete_ham}
\end{equation}
where $\Dtd \qk$ is defined implicitly as a function of $(\qdk,\pdk)$ through equation
(\ref{eq:m_legendre}). $\bar{H}_d$ is related to the previously defined Hamiltonian function by the
following relationship:
\begin{equation}
\bar{H}_d(\qdk,\pdk)=\Dtd \tdk H_d(\qdk,\pdk)\,.
\end{equation}
In addition, we have: $\edk=-H_d(\qdk,\pdk)$, that is, the momentum associated with the time is the
opposite of the Hamiltonian.

 Let $S_d^H$ be the discrete action summation:
\begin{eqnarray}
S^H_d&=&\tau\sum_{k=0}^{n-1} \langle\pdk, \Dtd\qdk \rangle-\bar{H}_d(\qdk,\pdk) \\
&=&\tau\sum_{k=0}^{n-1} \langle\pdk, \Dtd\qdk \rangle+\langle \edk, \Dtd \tdk \rangle\,.
\end{eqnarray}

Before stating the generalized discrete modified Hamilton's principle, we need to remark that all
the coordinates are not independent since the \emph{holonomic} constraint $\edk=-H(\qdk,\pdk)$ holds. There are two ways to handle this situation \cite{blo03bail}, one can either replace $\edk$ by $-H(\qdk,\pdk)$ in the action and then take the variations or one can use Lagrange multiplier to append the constraint $\edk+H(\qdk,\pdk)=0$ to the integral. Therefore we can give two equivalent formulations of the GDMHP.
\begin{defn}[Generalized discrete modified Hamilton's principle]
Critical points of the discrete action
$$S_d^H=\tau\sum_{k=0}^{n-1} \langle\pdk, \Dtd\qdk \rangle+\langle \edk, \Dtd \tdk \rangle$$ 
in the class of discrete curves $(\qdk,\tdk,\pdk,\edk)$ with endpoints $(\tau_0,t_0, q_0)$ and $(\tau_n,t_n,q_n)$ subject to the constraint $\edk+H_d(\qdk,\pdk)=0$ correspond to discrete trajectories of the discrete Hamiltonian system going from $(t_0,q_0)$ to $(t_n,q_n)$.
\label{def:m_discrete_pcp_la_v2}
\end{defn}
\begin{defn}[Generalized discrete modified Hamilton's principle] $\,$ $\:$ $\,$
Critical points of the discrete action 
$$S_d^H=\tau \sum_{k=0}^{n-1} \langle\pdk, \Dtd\qdk \rangle-H_d(\qdk,\pdk)\Dtd \tdk$$
in the class of discrete curves $(\qdk,\tdk,\pdk)$ with endpoints $(\tau_0,t_0, q_0)$ and $(\tau_n,t_n,q_n)$ correspond to trajectories of the discrete Hamiltonian system going from $(t_0,q_0)$ to $(t_n,q_n)$.
\label{def:m_discrete_pcp_la}
\end{defn}
To derive the equations of motion we need to specify the discrete derivative operator, $\Dtd$ and
its associated Leibnitz law.

\subsection{Examples}
\subsubsection{St{\"o}rmer type of algorithm}
\paragraph{Lagrangian approach}
Consider a Lagrangian function $L(q,\dot{q})$ and define the discrete Lagrangian map trivially by $L_d(\qk,\Dt
\qk)=L(\qk, \Dt \qk)$. Discrete equations of motion are obtained from the generalized discrete
Hamilton's principle:
\begin{eqnarray}
\delta S_d^L&=&\tau\sum_{k=0}^{n-1} \delta \bar{L}_d(\qk,\tk, \Dt \qk,\Dt\tk)\nonumber\\
&=&\tau\sum_{k=0}^{n-1} \delta (\Dt \tk L_d(\qk,\frac{\Dt \qk}{\Dt \tk}))\nonumber\\
&=&\tau\sum_{k=0}^{n-1} (\delta \Dt \tk ) L_d^k + \Dt \tk D_1 L_d^k \delta \qk\nonumber\\
&&+\Dt\tk D_2L_d^k
 \left(\frac{\Dt \delta \qk}{\Dt \tk}-\frac{\Dt\qk}{(\Dt
\tk)^2}\delta \Dt \tk\right)\nonumber\,,
\end{eqnarray}
where  $L_d^k=L_d(\qk,\frac{\Dt
\qk}{\Dt \tk})$. Using the Leibnitz law (equation \eqref{eq:modified_leibnitz}) and the fixed end points constraint, we obtain:
\begin{equation}
\delta S_d^L=\tau\sum_{k=1}^{n-1} -\Dt e_{k-1}\delta \tk+(-\Dt D_2L_d^{k-1}+\Dt \tk D_1 L_d^k )\delta \qk\,,
\end{equation}

where we have used the fixed end points constraint to derive the last equation and defined
$$\eko=\pfrac{\bar{L}_d^k}{\Dt \tk}=L_d(\qk,\frac{\Dt \qk}{\Dt
\tk})-D_2L_d (\qk,\frac{\Dt \qk}{\Dt \tk})\frac{\Dt\qk}{\Dt\tk}\,.$$

Finally we obtain the modified Euler-Lagrange equations by setting the variations to zero:
\begin{eqnarray}
\ek-e_{k-1}&=&0\nonumber\,,\\
\Dt \tk D_1 L_d(\qk,\frac{\Dt
\qk}{\Dt \tk})-\Dt D_2L_d(\qkmo,\frac{\Dt
\qkmo}{\Dt t_{k-1}})&=&0\,.
\label{eq:m_l_stormer}
\end{eqnarray}

\begin{lem}
The algorithm defined by (\ref{eq:m_l_stormer}) preserves the Lagrangian two-form and the energy.
\end{lem}
\begin{pf}
The first equation of the algorithm proves energy conservation. To show that the Lagrangian
two-form is preserved, we compute $dS^L_d$ along a discrete trajectory:
\begin{eqnarray}
dS^L_d&=&
\tau\sum_{k=1}^{n-1} \Dt (L_d^{k-1} d \tk)+
\Dt(D_2L_d^{k-1}d \qk)-\Dt (\frac{D_2L_d^{k-1}}{\Dt t_{k-1}} \Dt\qkmo d \tk)\nonumber\\
&=&\tau\sum_{k=1}^{n-1} \Dt (\ek d\tk+D_2L_d^{k-1}d \qk) \nonumber\\
&=&\tau\sum_{k=1}^{n-1} \Dt \theta_k^L\,,
\end{eqnarray}
where $\theta_k^L=\ek d\tk+D_2L_d^{k-1}d \qk$. Since $d^2=0$, we obtain that the symplectic
two-form $\omega_k^L=d\theta_k^L$ is preserved along the trajectory.
\end{pf}

The proof of this lemma only involves the modified Leibnitz law and does not depend on the
definition of the discrete Lagrangian function. As a consequence, it also applies if one derives
modified velocity Verlet and Newmark algorithms.

\paragraph{Hamiltonian approach}
Let the Lagrangian function be $L(q,\dot{q})=\undemi \dot{q}^T M \dot{q}-V(q)$. Then
\begin{equation}
\bar{L}_d=\Dt \tk \left(\undemi \frac{\Dt \qk}{\Dt \tk} M \frac{\Dt \qk}{\Dt \tk}-V(\qk)\right)\,,
\end{equation}
and the associated momenta are:
\begin{eqnarray}
\pko&=&M  \frac{\Dt \qk}{\Dt \tk}   \nonumber\,,\\
e_{k+1}&=&-\undemi \frac{\Dt \qk}{\Dt \tk} M \frac{\Dt \qk}{\Dt \tk}-V(\qk)\,.
\end{eqnarray}
The discrete Hamiltonian function is then:
\begin{equation}
\bar{H}_d=\Dt \tk (\undemi \pko^T M^{-1} \pko +V(\qk))=\Dt t_k H_d(\qk,\pko)\,.
\end{equation}
One can readily verify that $H_d(\qk,\pko)=-e_{k+1}$.

Let us now derive the modified discrete equations of motion by applying the GDMHP (theorem
(\ref{def:m_discrete_pcp_la})). We skip a few steps in the evaluation of the variations of $S_d^H$
to finally find:
\begin{eqnarray}
\delta S_d^H&=&\tau\delta \sum_{k=0}^{n-1} \langle \pko,\Dt \qk\rangle -\bar{H}_d(\qk,\pko)\nonumber\\
&=&\tau\sum_{k=0}^{n-1} \langle \Dt \qk- \Dt \tk D_2 H_d(\qk,\pko),\delta \pko \rangle\nonumber\\
&& -\langle \Dt \pk
+\Dt\tk D_1H_d(\qk,\pko),\delta \qk\rangle+\Dt \eko \delta \tko\,.\nonumber\\
\end{eqnarray}
The variations of $(\delta\qk,\delta\pko,\delta\tk)$ being independent, we obtain:
\begin{eqnarray}
\Dt \qk&=& \Dt\tk \pko \nonumber\,,\\
\Dt \pk&=&-\Dt\tk \nabla V(\qk)\nonumber\,,\\
\Dt \ek&=&0\,.
\label{eq:m_h_stormer}
\end{eqnarray}
\begin{lem}
The algorithm defined by equations (\ref{eq:m_h_stormer}) preserves the symplectic two-form and the
energy.
\end{lem}
\begin{pf}
The proof proceeds as the previous ones, we compute $dS_d^H$ along a discrete trajectory. We skip
the detail of the computation:
\begin{equation}
dS^H_d=\tau\sum_{k=0}^{n-1} \Dt \langle \pk,d\qk \rangle+\ek d \tk\,.
\end{equation}
Define $\theta_k^H= \langle \pk,d\qk \rangle+\ek d \tk\ $ and $\omega^H_k=d\theta_k^H$. Since
$d^2=0$, we obtain $$\Dt \omega_k^H=0\,.$$
\end{pf}
\begin{rem}
The $1$-form $\theta^H_k$  corresponds to the contact $1$-form $\theta$ encountered in continuous
time dynamics. Indeed, if one remembers that $\ek=-H_d(\qkmo,\pk)$, then we have:
\begin{eqnarray}
\theta &=&p dq -Hdt\,,\\
\theta_k^H&=& \pk d \qk -H_d(\qkmo,\pk)d\tk\,.
\end{eqnarray}
\end{rem}

\subsubsection{Midpoint discretization}
In the same manner, we can apply the modified variational principle to other discretization. For
the midpoint scheme we have $\qdk=\frac{\qko+\qk}{2}$  and the modified Leibnitz rule is defined by
equation (\ref{eq:leibnitz_midpoint}). Let us define the generalized momenta:
\begin{eqnarray}
\frac{\pko+\pk}{2}=&\pdk&=\pfrac{\bar{L}_d}{\Dtd \qdk}\,,\\
\frac{\eko+\ek}{2}=&\edk&=\pfrac{\bar{L}_d}{\Dtd \tdk}\,.
\end{eqnarray}
Then applying the modified discrete Hamilton's principle (Definition (\ref{def:m_discrete_pcp_la}))
yields (after a few simplifications):
\begin{equation}
\delta S_d^H=\tau\sum_{k=0}^{n-1} \langle \Dtd \tdk D_1 L_d^k -\Dtd \pdk,\delta \qdk \rangle -\Dtd
\edk \delta \tdk\,,
\end{equation}
where $L_d^k=L_d(\qdk,\frac{\Dtd \qdk}{\Dtd \tdk})$. The variations $(\delta \qdk,\delta \tdk)$
being independent, we obtain:
\begin{eqnarray}
\frac{\pko-\pk}{\tau}&=&\frac{\tko-\tk}{\tau}D_1L_d(\frac{\qko+\qk}{2},\frac{\qko-\qk}{\tko-\tk})\nonumber\,,\\
\eko&=&\ek\nonumber\,,\\
\frac{\pko+\pk}{2}&=&\frac{\tko-\tk}{\tau}D_2L_d(\frac{\qko+\qk}{2},\frac{\qko-\qk}{\tko-\tk})\nonumber\,,\\
\frac{\eko+\ek}{2}&=&L_d(\frac{\qko+\qk}{2},\frac{\qko-\qk}{\tko-\tk})\nonumber\\
&&-\langle D_2L_d(\frac{\qko+\qk}{2},\frac{\qko-\qk}{\tko-\tk}),\frac{\qko-\qk}{\tko-\tk}\rangle\,.
\label{eq:m_l_midpoint}
\end{eqnarray}

\begin{lem}
The algorithm defined by equations (\ref{eq:m_l_midpoint}) preserves the Lagrangian two-form as
well as the energy.
\end{lem}
\begin{pf}
We omit the proof since it proceeds as before.
\end{pf}

Now define the discrete Hamiltonian function
$H_d(\qdk,\pdk)=H(\frac{\qko+\qk}{2},\frac{\pko+\pk}{2})$ and the modified Hamiltonian function
$\bar{H}_d=\Dtd \tdk H_d(\qdk,\pdk)$. Then applying the GDMHP yields:
\begin{eqnarray}
\frac{\qko-\qk}{\tau}&=&\frac{\tko-\tk}{\tau} D_2 H_d(\frac{\qko+\qk}{2},\frac{\pko+\pk}{2})\,,\nonumber\\
\frac{\pko-\pk}{\tau}&=&-\frac{\tko-\tk}{\tau} D_1 H_d(\frac{\qko+\qk}{2},\frac{\pko+\pk}{2})\,,\nonumber\\
\eko-\ek&=&0\nonumber\,,\\
\frac{\eko+\ek}{2}&=&-H_d(\frac{\qko+\qk}{2},\frac{\pko+\pk}{2})\,.
\label{eq:m_h_midpoint}
\end{eqnarray}

\begin{lem}
The algorithm defined by equations (\ref{eq:m_h_midpoint}) preserves the symplectic two-form as
well as the energy.
\end{lem}
\begin{pf}
We omit the proof since it proceeds as before.
\end{pf}

\subsection{Concluding remarks}

The algorithm defined by equations (\ref{eq:m_h_midpoint}) is the same as the one developed by
Shibberu \cite{shi92}. Shibberu's approach corresponds to the first formulation of the GDMHP (definition
(\ref{def:m_discrete_pcp_la_v2})) for the midpoint rule but he used a different discrete variational principle from DVPII.

One other work on symplectic energy preserving algorithms is that of Kane, Marsden and Ortiz
\cite{kan99}. They developed a generalized discrete modified Hamilton's principle that is based on
DVPI. Their approach is different from ours: they assume a different time step at each iteration, and
then take the variation of the discrete action without varying the time step (i.e., in a $n$
dimensional space). As a consequence they only obtain $n$ equations for the $n+1$ variables
$(\qk,h_k)$ where $h_k$ is the time step at the $k^{th}$ step. They then add an energy constraint
to obtain $n+1$ equations. Their definition of the energy is similar to ours and therefore both
methods provide the same algorithms. However, there are fundamental differences between the two
methods. First, the method developed in this paper is fully variational.
Second, all the differences between DVPI and DVPII that we emphasize at the beginning of this paper
still remain because their work is based on DVPI whereas our is based on DVPII.

\section{Discrete Hamilton-Jacobi theory}
\label{sec:dhj}
So far we have developed two variational principles that are the discrete counterparts of
Hamilton's principle on the tangent bundle and on the cotangent bundle. Through several examples we
have observed that both variational principles are equivalent and that they allow us to recover
classical variational symplectic integrators. We have also shown that they can be modified so that
energy conservation is assured. In this section, we concentrate on discrete Hamilton-Jacobi theory.
We define discrete canonical transformations (DCT), discrete generating functions (DGF) and derive
a discrete Hamilton-Jacobi equation that allows us to show that the energy error for a certain class of scheme is invariant under discrete canonical transformations.

\subsection{Discrete symplectic geometry}
We consider again a discretization of the time $t$ into $n$ instants ${\mathcal T}=\{(t_k)_{k\in
[1,n]}\}$ but we restrict here to the case where $M_k$ is a $n$-dimensional vector space. We still
define ${\mathcal M}=\bigcup M_k$.

\begin{defn}
A discrete symplectic form $\omega$ on $\mathcal M$ is such that at $t_k$, $\omega=\omega_k^d$,
where $\omega_k^d$ is a non degenerate, closed, two-form on $M_k^d=M_k\cup M_{k+1}$.
\newline A discrete canonical one-form, $\theta$ on $\mathcal M$ is such that at $t_k$, $\theta=\theta_k^d$,
and $\omega_k^d=-d\theta_k^d$.
\newline A discrete symplectic vector space $({\mathcal M},\omega)$ is a vector space ${\mathcal M}=\bigcup
M_k$ together with a discrete symplectic two form on $\M$.
\end{defn}

Using a symplectic chart, a discrete symplectic form on $\mathcal M$ at $t_k$ can be written as:
\begin{equation}
\omega_k^d=d\qdk\wedge d\pdk\,,
\end{equation}
and the canonical one-form as $\theta_k^d=\pdk d\qdk$.

In the remainder of this section we consider the geometry associated with the midpoint scheme, that
is, we define $z_k^d=(\qdk,\pdk)$ as $z^d_k=\frac{\zk+\zko}{2}$ and use the modified Leibnitz law
\eqref{eq:leibnitz_midpoint}. However, the content of this section can be applied to any scheme as
long as one can define a discrete Hamiltonian vector field from the discrete Hamiltonian function
and the discrete symplectic two-form (see next definition). It is clear that the theory herein can
be adapted to systems for which the action integral involves a term of the form $H_d(z^d_k)$, where
$z^d_k$ is a linear combination of $z_k$ and $z_{k+1}$ but it is not clear if it can be adapted to
the St\"ormer rule for instance ($\zdk=(\qk,\pko)$ cannot be written as a linear combination of
$\zko$ and $\zk$ so the next definition does not apply). We do not know how to
modify this approach so that a discrete Hamiltonian vector field can be defined from the Hamiltonian function
$H_d(\qk,\pko)$.

\begin{defn}
Let $({\mathcal M},\omega)$ be a discrete symplectic vector space, and $H_d : \mathcal M\rightarrow
\mathbb{R}$ a smooth function. Define the discrete vector field $X_H^d$ such that at $t_k$, $X_H^d=X_k^d$,
where $X_k^d$ is of the form
\begin{equation}
X_H^d=\Dtd \qdk \pfrac{}{\qdk}+\Dtd \pdk \pfrac{}{\pdk}\,,
\end{equation}
and verifies:
\begin{equation}
i_{X_k^d}\omega_k^d=dH_d\,.
\label{eq:vect_field}
\end{equation}
The discrete vector field $X_H^d$ is called the discrete Hamiltonian vector field. \newline
$({\mathcal M},\omega,X_H^d)$ is called a discrete Hamiltonian system.
\end{defn}

\begin{prop}
Using the canonical coordinates, a Hamiltonian vector field is of the form:
\begin{equation}
X_H^d=J\cdot dH_d\,,
\label{eq:vect_field_coord}
\end{equation}
\end{prop}
where $J=\begin{pmatrix} 0 & I\\-I&0\end{pmatrix}$.
\begin{pf}
Equation (\ref{eq:vect_field}) is expressed in local coordinates as:
\begin{equation}
i_{X_H^d}(d\qdk\wedge d\pdk)=D_1H_d(\zdk)d\qdk +D_2H_d(\zdk)d\pdk\,.
\label{eq:vect_field_coord2}
\end{equation}
Let $X_H^d$ be:
\begin{equation}
X_H^d=\Dtd \qdk \pfrac{}{\qdk}+\Dtd \pdk \pfrac{}{\pdk}\,,
\end{equation}
then,
\begin{eqnarray}
i_{X_H^d}(d\qdk\wedge d\pdk)&=&(i_{X_H^d}d\qdk)d \pdk -d\qdk \wedge (i_{X_H^d}d \pdk)\\
&=&\Dtd \qdk d\pdk -\Dtd \pdk d\qdk\,.
\end{eqnarray}
Identifying this last equation with equation (\ref{eq:vect_field_coord2}) leads to equation
(\ref{eq:vect_field_coord}).
\end{pf}

\subsection{Discrete canonical transformation}
We now define the class of discrete canonical transformations. The definition given here is restricted
to linear time-dependent maps (with respect to the phase space variables). We believe larger class
of transformations may be considered if one works with discretization of the spacetime \cite{mar98pat}. Let $({\mathcal M}_1,\omega_1)$ and $({\mathcal M}_2,\omega_2)$ be discrete symplectic vector spaces and $\mathcal F$ be the set maps $f: {\mathcal T}\times {\mathcal M}_1 \rightarrow {\mathcal T}\times {\mathcal M}_2$ that are linear with respect to the phase space variables. Consider a map $f\in \mathcal F$ such that
$\forall t_k\in {\mathcal T}$, $f(t_k,\cdot)=f_k(\cdot)$ where $f_k$ is the following linear map:
\begin{eqnarray}
M_{1,k}^d& \rightarrow &  M_{2,k}^d\nonumber\\
\zk=(\qk,\pk)& \mapsto & Z_k=(\Qk,\Pk)=A_k \zk+B_k \nonumber\,.
\end{eqnarray}
Since $f_k$ is linear, we have:
\begin{eqnarray}
f_k(\zdk)&=&\undemi (f_k(z_k)+f_k(z_{k+1}))\label{eq:dct_zdk}\,,\\
f_k(\Dtd \zdk)&=& A_k \Dtd \zdk\label{eq:dct_dtd}\,.
\end{eqnarray}

\begin{defn}
A linear, time-dependent map $f$ is called a discrete canonical transformation (DCT) (or a discrete symplectic map) if
and only if $f^*\omega_2=\omega_1$, or equivalently, $\forall k\in [1,n]$,
$f_k^*\omega_{2,k}^d=\omega_{1,k}^d$.
\end{defn}

\begin{prop}
If $f$ is a DCT then $A_k$ is invertible for all $k\in [1,n]$
\end{prop}
\begin{pf}
Suppose there exists a $k$ such that $A_k$ is not invertible. Then $\exists \zdk \in M_{1,k}^d$
such that $$\exists v_1\in T_{\zdk}M_{1,k}^d |Tf_k\cdot v_1=0\,.$$ Then, $\forall v_2\in
T_{\zdk}M_{1,k}^d | v_2\ne 0$, $\omega_{1,k}^d(v_1,v_2)=\omega_{2,k}^d(Tf_k\cdot v_1,Tf_k\cdot
v_2)$ since $f$ is symplectic. The right hand side is zero but the left hand side is not. This is a
contradiction and therefore $A_k$ is invertible.
\end{pf}

\begin{lem}
Let $f$ be a discrete canonical transformation. Then $f_k^*\omega_{2,k}^d=\omega_{1,k}^d$ can be
written in the matrix form $A_kJA_k^T=J$. In addition, $f$ preserves the form of the discrete
Hamilton's equations.
\end{lem}
\begin{pf}
The statement $A_kJA_k^T=J$ is just the matrix statement of $f_k^*\omega_{2,k}^d=\omega_{1,k}^d$.
Let us prove that $f$ preserves the form of the discrete Hamilton's equations. Define the function
$K_d$ such that $K_d\circ f=H_d$.

On one hand, using equation (\ref{eq:dct_dtd}) we have:
\begin{eqnarray}
\Dtd Z_k^d&=&\frac{f_k (\zko)- f_k (\zk)}{\tau}\\
&=& A_k \Dtd \zdk\,.
\end{eqnarray}
On the other hand:
\begin{eqnarray}
J\nabla H_d(\zdk)&=&J \nabla (K_d \circ f_k(\zdk))\\
&=&J A_k^T \nabla K_d(\zdk)\,.
\end{eqnarray}
Since  $A_k J A_k^T=J$, we obtain:
\begin{equation}
\Dtd Z_k^d=J\nabla K_d(\zdk)
\end{equation}
\end{pf}
This last result can be summarized as follows:
\begin{prop}
Let $X_H^d$ be a discrete Hamiltonian vector field with Hamiltonian function $H_d$ and $f$ a
discrete symplectic map. Then $f_* X_H^d $ is a discrete Hamiltonian vector field with Hamiltonian
function $f_*H_d$.
\label{prop:sum}
\end{prop}

\subsection{Discrete generating functions}

\begin{prop}
Let $({\mathcal M}_1,\omega_1)$ and $({\mathcal M}_2,\omega_2)$ be two discrete symplectic vector
spaces, $\pi_{i}: {\mathcal M}_1\times {\mathcal M}_2\rightarrow {\mathcal M}_i $ the projection
onto ${\mathcal M}_i$  and define
\begin{equation}
\Omega=\pi_1^*\omega_1-\pi_2^*\omega_2\,.
\end{equation}
Then,
\begin{enumerate}
\item $\Omega$ is a discrete symplectic form on ${\mathcal M}_1\times {\mathcal M}_2$,
\item a map $f:{\mathcal M}_1\rightarrow {\mathcal M}_2$ is a discrete symplectic map if and only if $i_f^*\Omega=0$, where
$i_f: \Gamma_f\rightarrow {\mathcal M}_1\times {\mathcal M}_2$ is the inclusion map and $\Gamma_f$
is the graph of $f$.
\end{enumerate}
\end{prop}
\begin{pf}
We recall that at $t_k$, $\Omega=\Omega_k^d$ where
$\Omega_k^d=\pi_{1k}^*\omega_{1k}^d-\pi_{2k}^*\omega_{2k}^d$. To prove that $\Omega $ is a discrete
symplectic form, we need to prove that $\Omega_k^d$ is a symplectic form on $M_{1,k}^d \times
M_{2,k}^d$ for all $k\in [1,n]$.
\begin{eqnarray}
d\Omega_k^d&=&d(\pi_1^*\omega_{1,k}^d-\pi_2^*\omega_{2,k}^d)\\
&=&\pi_1^*d\omega_{1,k}^d -\pi_2^*d\omega_{2,k}^d\\
&=&0\,,
\end{eqnarray}
since $\omega_{i,k}^d$ is closed and $d$ commutes with the pull back operator.

Now let $\zdk=(z_{1,k}^d,z_{2,k}^d)\in M_{1,k}^d \times M_{2,k}^d$ and $v=(v_1,v_2)\in
T_{\zdk}(M_{1,k}^d \times M_{2,k}^d)\sim T_{z_{1,k}^d}M_{1,k}^d \times T_{z_{2,k}^d}M_{2,k}^d$ such
that
\begin{equation}
\forall w=(w_1,w_2)\in T_{\zdk}(M_{1,k}^d \times M_{2,k}^d)\quad \Omega_k^d(v,w)=0
\end{equation}
and let us prove that $v$ is zero. We have
\begin{eqnarray}
\Omega_k^d(v,w)&=&\omega_{1,k}^d (\pi_1 (\zdk))(T\pi_1\cdot v,T\pi_1\cdot w)-
\omega_{2,k}^d (\pi_2 (\zdk))(T\pi_2\cdot v,T\pi_2\cdot w)\\
&=&\omega_{1,k}^d (z_{1,k}^d)( v_1, w_1)-\omega_{2,k}^d (z_{2,k}^d)( v_2, w_2)
\label{eq:a1}
\end{eqnarray}
The  right hand side of equation (\ref{eq:a1}) is zero for all $w$ if and only if both terms are
zero, that is,
\begin{equation}
\omega_{1,k}^d (z_{1,k}^d)( v_1, w_1)=0\,,\;
\omega_{2,k}^d (z_{2,k}^d)( v_2, w_2)=0
\end{equation}
Since $\omega_{i,k}^d$ is non degenerate, $v_1=v_2=0$ and $\Omega_k^d$ is closed.

We now prove the second statement of the proposition. We first notice that $f_k$ induces a
diffeomorphism of $M_{1,k}^d$ to $\Gamma_{f_k}$, so we can write
\begin{equation}
T_{(\zdk,f_k(\zdk))}=\left\{(v,Tf_k\cdot v)|v\in T_{\zdk}M_{1,k}^d\right\}
\end{equation}
Then,
\begin{eqnarray}
i^*\Omega_k^d((v_1,Tf_k\cdot v_1),(v_2,Tf_k\cdot v_2))&=&
\omega_{1,k}^d(v_1,v_2)-\omega_{1,k}^d(Tf_k\cdot v_1,Tf_k\cdot v_2)\nonumber\\
&=&(\omega_{1,k}^d-f_k^*\omega_{2,k}^d)(v_1,v_2)
\end{eqnarray}
Hence, $f_k$ is symplectic if and only if $i^*\Omega_k^d=0$, i.e., $f$ is a discrete symplectic map
if and only if $i^*\Omega =0$.
\end{pf}

Using the Poincar{\'e} lemma we may write $\Omega_k^d=-d\Theta_k^d$ and the previous proposition
says that $i_{f_k}^*\Theta_k^d$ is closed if and only if $f$ is a discrete symplectic map. Using
again the Poincar{\'e} lemma, we conclude that if $f$ is a discrete symplectic map then there
exists a function $S:\Gamma_f\rightarrow\mathbb{R}$ such that $i_f^*\Theta=dS$, i.e.,  $\forall
k\in [1,n]$, $i_{f_k}^*\Theta_k^d=dS_k$

\begin{defn}
Such a function $S$ is called a discrete generating function for the discrete symplectic map $f$.
$S$ is locally defined and depends on the choice of $\Theta$.
\end{defn}


\begin{itemize}
\item Let $\theta_{1,k}^d=\pdk d\qdk$ and $\theta_{2,k}^d=\Pdk d\Qdk$, then
\begin{eqnarray}
i^*_{f_k}\Theta_k^d&=&(\pi_1\circ i_{f_k})^*\pdk d\qdk-(\pi_2\circ i_{f_k})^*\Pdk d\Qdk\,, \\
d S&=&\pfrac{S}{q}(\qdk,\Qdk)d\qdk+\pfrac{S}{Q}(\qdk,\Qdk)d\Qdk\,,
\end{eqnarray}
that is,
\begin{equation}
\pdk=\pfrac{S}{q}(\qdk,\Qdk) \qquad \Pdk =-\pfrac{S}{Q}(\qdk,\Qdk)\,.
\end{equation}
$S$ as defined corresponds to a discrete generating function of the first kind.

\item Let $\theta_{1,k}^d=\pdk d\qdk$ and $\theta_{2,k}^d=-\Qdk d\Pdk$, then
\begin{eqnarray}
i^*_{f_k}\Theta_k^d&=&(\pi_1\circ i_{f_k})^*\pdk d\qdk+(\pi_2\circ i_{f_k})^*\Qdk d\Pdk\,, \\
d S&=&\pfrac{S}{q}(\qdk,\Qdk)d\qdk+\pfrac{S}{Q}(\qdk,\Qdk)d\Qdk\,,
\end{eqnarray}
that is,
\begin{equation}
\pdk=\pfrac{S}{q}(\qdk,\Qdk) \qquad \Qdk =\pfrac{S}{P}(\qdk,\Pdk)\,.
\end{equation}
$S$ as defined corresponds to a discrete generating function of the second kind.
\end{itemize}
In the same way, one can define $4^n$ generating functions as in the continuous case. Note that
since $f$ is linear with respect to its spatial variables, $S$ is also linear with respect to its
spatial variables. At $t_k$, $S=S_k$ where $S_k(\cdot)=T_k (\cdot)+U_k$
is affine map, $T_k$ is a $ 2n\times 2n $
matrix and $U_k$ is a $2n\times 1 $ matrix.

\subsection{Discrete Hamilton-Jacobi theory}

In this section we use the notions introduced previously to develop a discrete Hamilton-Jacobi
theory. Let $f$ be a discrete symplectic map, let
 $M_{i,k}^d=T^*{\mathcal Q}_{i,k}^d$ and let $S$ be an
associated discrete generating function such that at $t_k$, $S=S_k^d$ where $S_k(\cdot)=T_k (\cdot)+U_k$

\begin{thm}
Define $$\tpdk(\qdk,\Qdk)=D_1S_k(\qdk,\Qdk)\,,\quad \tPdk(\qdk,\Qdk)=-D_2S_k(\qdk,\Qdk)\,.$$ Then
the following two conditions are equivalent:
\begin{enumerate}
\item $S$ is a discrete generating function associated with $f$;
\item \begin{itemize}\item For every curve $(c_k)_k$ in ${\mathcal Q}_1=\bigcup {\mathcal Q}_{1,k}$ satisfying:
\begin{equation}
\Dtd \cdk= T\pi_{{\mathcal Q}_{1,k}^d}^* X_H^d(\cdk,\tpdk)\,,
\end{equation}
the curve $k \mapsto (\cdk,\tpdk)$ is a discrete integral curve of $X_H^d$, where $\pi_{{\mathcal
Q}_{1,k}^d}^*$ is the cotangent bundle projection onto the configuration space.
\item For every curve $(c_k)_k$ in ${\mathcal Q}_2=\bigcup {\mathcal Q}_{2,k}$
satisfying:
\begin{equation}
\Dtd \cdk= T\pi^*_{{\mathcal Q}_{2,k}^d} X_K^d(\cdk,\tPdk)\,,
\end{equation}
the curve $k \mapsto (\cdk,\tPdk)$ is a discrete integral curve of $X_K^d$, where $\pi^*_{{\mathcal
Q}_{2,k}^d}$ is the cotangent bundle projection onto the configuration space.
\end{itemize}
\end{enumerate}
\end{thm}

\begin{pf}
Suppose $S$ is a discrete generating function, let $\Qdk$ be fixed
and consider a curve
$(c_k)_k$ verifying
\begin{equation}
\Dtd \cdk= T\pi_{{\mathcal Q}_{1,k}^d}^* X_H^d(\cdk,\tpdk)\,,
\end{equation}
In other words, $\ck$ verifies:
\begin{equation}
\Dtd \cdk = D_2 H(\cdk,\tpdk)\,,
\end{equation}
Since $S$ is a generating function, $\tpdk$ is the momentum associated with $\cdk$ and
 verifies:
\begin{equation}
\Dtd \tpdk = -D_1 H(\cdk,\tpdk)\,.
\end{equation}
These last two equations are exactly a restatement of: $k \mapsto (\cdk,\tpdk)$ is a discrete
integral curve of $X_H^d$. To derive the second item we proceed in the same manner, but this time
$\qdk$ is fixed.

Now we suppose item (2) and we show that $S$ is a discrete generating function for $f$. The
statements $k
\mapsto (\cdk,\tpdk)$ is a discrete integral curve of $X_H^d$  and $k
\mapsto (\cdk,\tPdk)$ is a discrete integral curve of $X_K^d$  are equivalent to saying that $\tpdk$ and
$\tPdk$ are the momenta associated with the generalized coordinates, and therefore, $S$ is a
generating function for $f$.
\end{pf}

\begin{thm}\label{theo:dhj}
We consider again a  time dependent function $S$ which is
linear with respect to the spatial variables. Then
the following two statements are equivalent:
\begin{enumerate}
\item $S$ is a discrete generating function associated with $f$;
\item For every $H$ there is a function $K$ such that
\begin{equation}
H(\qdk,D_1S(\qdk,\Qdk))=K(\Qdk,D_2S(\qdk,\Qdk))
\end{equation}
\end{enumerate}
\end{thm}
\begin{pf}
Suppose $S$ is a discrete generating function. Then from the previous theorem, for every curve
$(c_k,C_k)$ in ${\mathcal Q}_1\times {\mathcal Q}_2$ satisfying $\Dtd \cdk= T\pi_{{\mathcal
Q}_{1,k}^d}^* X_H^d(\cdk,\tpdk)$ and $\Dtd \Cdk= T\pi_{{\mathcal Q}_{2,k}^d}^* X_K^d(\Cdk,\tPdk)$,
the curves $k
\mapsto (\cdk,\tpdk)$ and $k
\mapsto (\Cdk,\tPdk)$ are discrete integral curves of $X_H^d$ and $X_K^d$ respectively. Then, using
the symplectic identity (\cite{abr78} page $382$) that holds for any function $S$
$$\omega_{1,k}^d (T(D_1S \circ \pi_{{\mathcal Q}_{1,k}^d}^*)\cdot v,w)=
\omega_{1,k}^d (v,w-T(D_1S \circ \pi_{{\mathcal Q}_{1,k}^d}^*)\cdot w)$$
we get:
\begin{eqnarray}
\omega_{1,k}^d (T&(&D_1S \circ \pi_{{\mathcal Q}_{1,k}^d}^*)\cdot X_H^d (c_k, D_1 S_k),w)=\nonumber\\
\omega_{1,k}^d &(&X_H^d (c_k, D_1 S_k),w)-dH_d(c_k, D_1 S_k)\cdot TD_1S(c_k, D_1 S_k)\cdot w\label{eq:tmp7}\\
\omega_{2,k}^d (T&(&-D_2S \circ \pi_{{\mathcal Q}_{2,k}^d}^*)\cdot X_K^d (C_k,- D_2 S_k),w)=\nonumber\\
\omega_{2,k}^d &(&X_K^d (C_k, -D_2 S_k),w)-dK_d(C_k, -D_2 S_k)\cdot T-D_2S(C_k, -D_2 S_k)\cdot
w\nonumber\\\label{eq:tmp8}
\end{eqnarray}
In addition, since $\pdk=D_1 S(\cdk,\Cdk)$ and $\Pdk=-D_1 S(\cdk,\Cdk)$,
\begin{eqnarray}
\Dtd \pdk&=&TD_1 S(\cdk,\Cdk) \Dtd \cdk=T(D_1S \circ \pi_{{\mathcal Q}_{1,k}^d}^*)\cdot X_H^d (c_k, D_1
S_k)\\
\Dtd \Pdk&=&T(-D_2S \circ \pi_{{\mathcal Q}_{2,k}^d}^*)\cdot X_K^d (C_k,- D_2 S_k)\,.
\end{eqnarray}
$f$ being a discrete canonical map, $Tf_k(\Dtd \pdk)=\Dtd \Pdk$ so the left hand side of equation
\ref{eq:tmp8} is the image under f of the left hand side of (\ref{eq:tmp8}). Using proposition
(\ref{prop:sum}), we conclude that:
$$Tf_k \cdot dH_d(c_k, D_1 S_k)\cdot TD_1S(c_k, D_1 S_k)=-dK_d(C_k, -D_2 S_k)\cdot TD_2S(C_k, -D_2 S_k)\,,$$
which is equivalent to the discrete Hamilton-Jacobi equation.

The proof that $2.$ implies $1.$ follows from these arguments.
\end{pf}

\subsection{Applications of the discrete Hamilton-Jacobi theory}
\begin{figure*}[htb]
\begin{center}
\subfigure[Trajectory of in the $q-p$ plane]
{\includegraphics[width=0.45\linewidth]{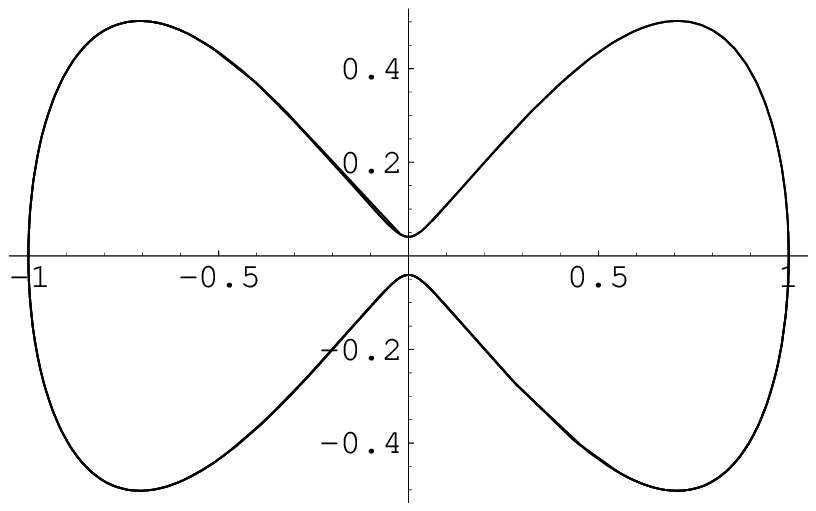}}
\subfigure[Energy error for constant time step midpoint scheme as a function of time.]
{\includegraphics[width=0.45\linewidth]{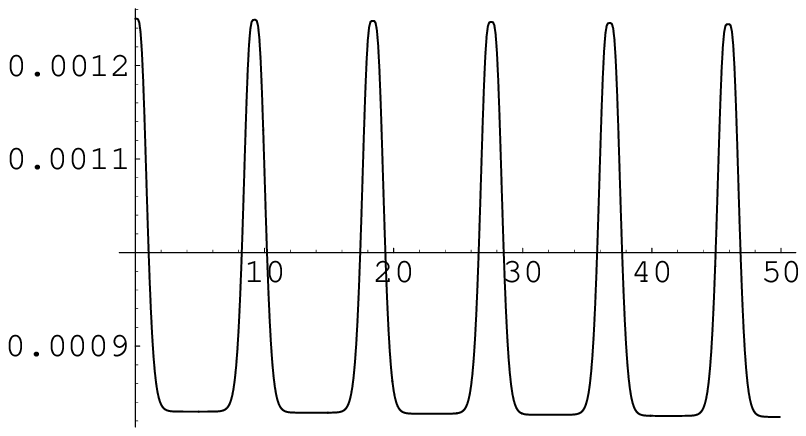}}
\caption{\label{fig:ex1_1}Particle in a double well potential with initial conditions $(q,p)=(1,0.05)$}
\end{center}
\end{figure*}

\begin{figure*}[htb]
\begin{center}
\subfigure[Trajectory of in the $q-p$ plane]
{\includegraphics[width=0.45\linewidth]{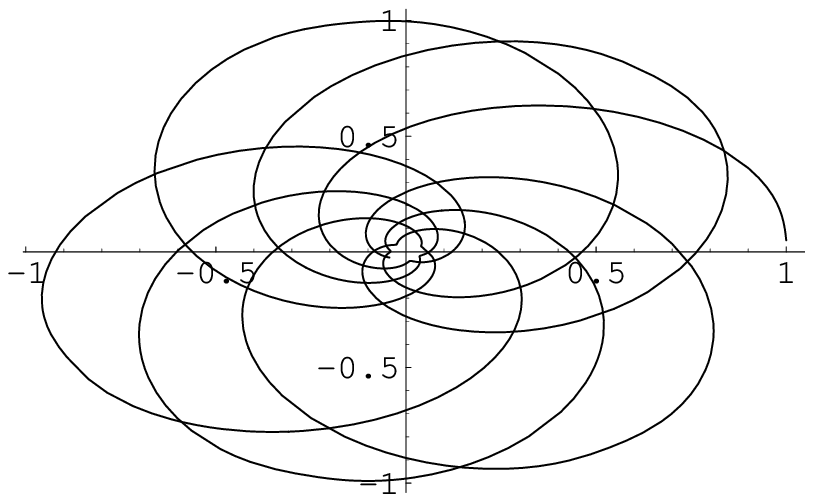}}
\subfigure[Energy error for constant time step midpoint scheme as a function of time.]
{\includegraphics[width=0.45\linewidth]{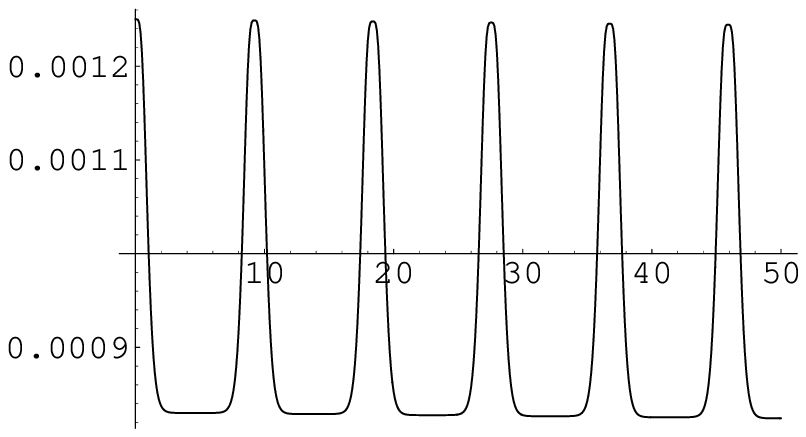}}
\caption{\label{fig:ex1_2}Particle in the Hamiltonian vector field $f_*X_H^d$, where $X_H^d$ is
the Hamiltonian vector field corresponding to a double well potential. Initial conditions are
$(Q,P)=f_0(1,0.05)$.}
\end{center}
\end{figure*}

The goal of this section is to highlight the benefit of having a discrete Hamilton-Jacobi theory.
First, we have proven the invariance of the discrete Hamilton's equations under a certain class of
coordinate transformations. Second, we have shown in theorem \ref{theo:dhj}
that changing coordinates using a
discrete symplectic map does not improve the performance of the algorithm in terms of energy
conservation. As a consequence we have the following lemma:
\begin{lem}
The midpoint scheme preserves the energy for linear systems.
\end{lem}
\begin{pf}
The discrete phase flow for linear systems is piecewise linear continuous and the map
$(\qk,\pk)\mapsto (\qko,\pko)$ is symplectic (the midpoint scheme is a symplectic algorithm).
Therefore, the discrete phase flow is a discrete symplectic map that maps $H$ into a constant $K$.
Integration of the new Hamiltonian system defined by $K$ is trivial ($(\Qko,\Pko)=(\Qk,\Pk)$) and
obviously preserves the energy. As a consequence the integration of the Hamiltonian system defined
by $H$ also preserves the energy.
\end{pf}
Finally, we illustrate the use of the above material with a nonlinear  example. We study the energy error in the integration of the equations of motion of a particle in a double well potential using different sets of canonical coordinates. Consider a particle in a
double well potential, i.e., $H=\undemi p^2+\undemi (q^4-q^2)$. As shown in figure
(\ref{fig:ex1_1}), the midpoint scheme does not preserve the energy.  The following time-dependent discrete canonical transformation (at each step the transformation is a different expression) $Z_k=A_kz_k+B_k$ where
$A_k=\begin{pmatrix} \cos(k \theta)& -\sin(k\theta)\\\sin(k \theta)& \cos(k\theta)\end{pmatrix}$, and $B_k=0$  rotates the system by $k\theta=k\arccos{0.99}$ at the $k^{th}$ step. In figure (\ref{fig:ex1_2}) we plot the same trajectory in the new system of coordinates, the energy error is exactly the same. In other words, the energy error is invariant under discrete canonical maps.

\section{Optimal control}
\label{sec:optimal}
For a general optimal control problem, necessary conditions for optimality may be derived from the
Pontryagin maximum principle. These conditions often yield equations of the same form as Hamilton's
equations coupled with nonlinear equations. We have seen previously that Hamiltonian systems, i.e.,
Hamilton's equations, can be integrated using symplectic integrators. However, if Hamilton's
equations are coupled with algebraic nonlinear equations, the above theory does not apply. What is
the correct discretization of the algebraic equation? In this section, we develop a discrete
maximum principle that tackles this problem and provides a unified view on solving optimal control
problems using symplectic integrators.

\subsection{Necessary conditions for optimality}
\subsubsection{Problem Statement}
Let $J=\int_0^T g(x,u)dt$ be a performance index (also called a cost function) and consider the
following optimal control problem:
\begin{eqnarray}
\min_u \int_{t_0}^{t_f} g(x,u)dt\,,
\label{eq:def_optimalpb}
\end{eqnarray}
subject to the dynamics
\begin{equation}
\dot{x}=f(x,u)\,,\label{eq:dyn}
\end{equation}
and to the initial and final time constraints:
\begin{equation}
\begin{array}{cc}
\phi_i(x(t_0),t_0)=0\,, &\phi_f(x(t_f),t_f)=0\,,\label{eq:boundcond}
\end{array}
\end{equation}
where $f$ and $g$ are functions from  $\mathbb{R}^n\times
\mathbb{R}^m$ to $\mathbb{R}$ of class $C^1$.

\subsection{Maximum principle}
To solve the optimal control problem, we apply the maximum principle.
\begin{thm}[Maximum principle]\label{theo:mp}
Solutions to the optimal control problem defined by equations \eqref{eq:def_optimalpb}, \eqref{eq:dyn} and \eqref{eq:boundcond} correspond to critical points of the cost function $J$ in the class of curves $\gamma=(x(t),u(t))\in \Gamma$
where $\Gamma$ is the set of curves satisfying (\ref{eq:dyn}) and (\ref{eq:boundcond}).
\end{thm}
\begin{rem}
This formulation differs from the one given by Pontryagin \cite{pon86} but the main point of the
Pontryagin maximum principle is that it yields necessary conditions for optimality under far less
severe regularity conditions. The above formulation is based on the equivalence between the
Pontryagin maximum principle and the calculus of variations in the case where the control region is
an open set in a finite dimensional vector space (see \cite{pon86} chapter V for more details). It is therefore equivalent to classical variational formulations given in Bloch et al. \cite{blo03bail,blo99cro} and Gregory and Lin \cite{gre92lin} for instance. 
\end{rem}
To apply the maximum principle we first need to define the augmented cost function $J_a$:
\begin{eqnarray}
J_a&=&\int_{t_0}^{t_f} g(x,u) + \langle p,\dot{x}-f(x,u)\rangle dt +\langle
\lambda_i ,\phi_i(x(t_0),t_0)\rangle+\langle \lambda_f ,
\phi_f(x(t_f),t_f)\rangle\nonumber\\
&=&\int_0^T  H(x,p,u)- \langle p,\dot{x}\rangle dt+\langle \lambda_i ,
\phi_i(x(t_0),t_0)\rangle+\langle \lambda_f ,
\phi_f(x(t_f),t_f)\rangle\,,\nonumber
\end{eqnarray}
where the $p$'s, the $\lambda_i$'s and the $\lambda_f$'s are Lagrange multipliers and
$H(x,p,u)=g(x,u)+\langle p,f(x,u)\rangle$. Taking variations of the augmented cost function
assuming fixed initial and final time yields:
\begin{eqnarray}
\delta J_a&=&\delta \left(\int_{t_0}^{t_f}  H(x,p,u) + \langle p,\dot{x}\rangle dt\right)+\delta\langle \lambda_i ,
\phi_i(x(t_0),t_0)\rangle\nonumber\\
&&+\delta \langle \lambda_f ,\phi_f(x(t_f),t_f)\rangle\nonumber\\
&=&\int_{t_0}^{t_f}\langle D_2H(x,p,u)-\dot{x},\delta p\rangle +\langle D_1H(x,p,u)+\dot{p},\delta
x\rangle\nonumber\\
&&{+}\:\langle D_3H(x,p,u),\delta u\rangle dt+\langle -p(t_f)+D_1\phi_f^T\lambda_f,\delta x_f\rangle
\nonumber\\
{+}\:\langle p(t_i)+D_1\phi_i^T\lambda_i,\delta x_i\rangle\,.\nonumber
\end{eqnarray}
We now let the variations of $J_a$ be zero to obtain necessary conditions for optimality:
\begin{eqnarray}
\dot{x}&=&D_2H(x,p,u)\,,\label{eq:neccond1}\\
\dot{p}&=&-D_1H(x,p,u)\,,\label{eq:neccond2}\\
0&=&D_3H(x,p,u)\label{eq:neccond3}\,,
\end{eqnarray}
as well as transversality conditions:
\begin{equation}
p(t_i)=-D_1\phi_i(x(t_0),t_0)^T\lambda_i\,,\;
p(t_f)= D_1\phi_f(x(t_f),t_f)^T\lambda_f\label{eq:trans}\,.
\end{equation}
Equations (\ref{eq:neccond1})-(\ref{eq:trans}) define the necessary conditions for optimality.

\subsection{Solving the necessary conditions for optimality}
To solve these conditions, the most common technique is to find the optimal control feedback law
from (\ref{eq:neccond3}) and then use a shooting method to solve the two-point boundary value
problem defined by (\ref{eq:neccond1}), (\ref{eq:neccond2}) and (\ref{eq:trans}). More precisely,
suppose (\ref{eq:neccond3}) allows one to solve for $u$ as a function of $(x,p)$ and define the
Hamiltonian function
\begin{equation}
\bar{H}(x,p)=H(x,p,u(x,p))\,,
\label{eq:barh}
\end{equation}
then the necessary conditions (\ref{eq:neccond1}) and (\ref{eq:neccond2}) simplify to:
\begin{eqnarray}
\dot{x}&=&D_2\bar{H}(x,p)\,,\label{eq:neccond1bis}\\
\dot{p}&=&-D_1\bar{H}(x,p)\,.\label{eq:neccond2bis}
\end{eqnarray}
Equations (\ref{eq:neccond1bis}) and (\ref{eq:neccond2bis}) define a Hamiltonian system that has
no physical meaning in general. As we will see later, for sub-Riemannian optimal control problems
the Legendre transform is ill-defined and therefore DVPI cannot be used to discretize such
systems whereas one could use DVPII (theorem \ref{def:discrete_pcp_la}). However, one may not be
able to solve (\ref{eq:neccond3}), and then the question of how one can use symplectic integrators
to solve the optimal control problem arises. What is the correct discretization of
(\ref{eq:neccond3})? In the next section we address this issue. Specifically, we introduce a discrete maximum principle that allows us to derive discrete necessary conditions for optimality that are in agreement with the one obtained from the maximum principle.

\subsection{Discrete maximum principle}
\subsubsection{Problem statement}

In discrete settings, the cost function is $$J=\sum_{k=0}^{n-1} g_d(\xdk,\udk)\tau\,,$$ and the
optimal control problem (\ref{eq:def_optimalpb}) is formulated as:
\begin{equation}
\min_{\udk} \sum_{k=0}^{n-1} g_d(\xdk,\udk)\tau\,,
\end{equation}
subject to the dynamics
\begin{equation}
\Dtd \xdk=f_d(\xdk,\udk)\,,
\label{eq:ddyn}
\end{equation}
and to boundary conditions:
\begin{equation}
\begin{array}{cc}
\phi_0(x_0,t_0)=0\,,&
\phi_n(x_n,t_n)=0\,,\label{eq:dboundcond}
\end{array}
\end{equation}
where $f_d$ and $g_d$ are functions from  $\mathbb{R}^n\times\mathbb{R}^m$ to $\mathbb{R}$ of class $C^1$. They correspond to discretization of the continuous time functions $f$ and $g$.

\subsubsection{Discrete maximum principle}
To obtain necessary conditions for optimality, we define the following discrete maximum principle, the discrete counterpart of the maximum principle:
\begin{defn}[Discrete maximum principle]\label{theo:dmp}
Solutions to the discrete optimal control problem correspond to critical points of the cost function $J$  in the class of discrete curves $\gamma\in\Gamma$, where $\Gamma$ is the set of all discrete curves $(\xk,\uk)_{k\in [1,n]}$ that verify (\ref{eq:ddyn}) and (\ref{eq:dboundcond}).
\end{defn}

\begin{rem}
The above definition is the discrete counterpart of the maximum principle. It compares to previous works on discrete optimal control theory that extend the Pontryagin maximum principle to discrete systems such as Jordan and Polak \cite{jor664pol} as theorem \ref{theo:mp} compares to the Pontryagin maximum principle. In other words, in contrast with Jordan and Polak \cite{jor664pol}, we restrict the class of discrete optimal control problems so that we can derive necessary conditions that define symplectic algorithms.
\end{rem}

As in the continuous case, to find critical points of $J$ under the non-holonomic constraint defined by equation \eqref{eq:ddyn}, we must append the constraints to $J$ using the Lagrange multipliers. The resulting function is called the augmented cost function:
\begin{eqnarray}
J_a&=&\sum_{k=0}^{n-1} (g_d(\xdk,\udk)-\langle{\pdk}, \Dtd \xdk -f_d(\xdk,\udk)\rangle)\tau
+\langle \lambda_0 ,\phi_0 \rangle+\langle \lambda_n,\phi_n\rangle\\
&=&\sum_{k=0}^{n-1} (H_d(\xdk,\pdk,\udk)-\langle{\pdk}, \Dtd \xdk\rangle)\tau+\langle
\lambda_0 ,\phi_0 \rangle+\langle \lambda_n,\phi_n\rangle\,,
\label{eq:Ja}
\end{eqnarray}
where the $\pk$'s, the $\lambda_0$'s and the $\lambda_n$'s are Lagrange multipliers and
$H_d(\xdk,\pdk,\udk)=g_d(\xdk,\udk)+\langle{\pdk},f_d(\xdk,\udk)\rangle$. To apply the discrete
maximum principle, one needs to specify the discrete derivative operator as well as the expressions
of $\xdk$, $\udk$ and $\pdk$ as a function of $(\xko,\xk)$, $(\uko,\uk)$ and $(\pko,\pk)$
respectively.

\subsubsection{Examples}
\paragraph{St{\"o}rmer's rule}
If we choose $\Dtd$ to be the forward difference $\Dt$ and $(\xdk,\pdk,\udk)=(\xk,\pko,\uk)$  then
we recover the discrete maximum principle developed by Bloch, Crouch, Marsden and Ratiu
\cite{blo02cro}.
\begin{eqnarray}
\delta J_a&=&\delta \left(\sum_{k=0}^{n-1} (H_d(\xdk,\pdk,\udk)+\langle{\pdk}, \Dtd
\xdk\rangle)\tau\right) +\delta\langle \lambda_0 ,\phi_0 \rangle+\delta\langle
\lambda_n,\phi_n\rangle\nonumber\\
&=&\sum_{k=0}^{n-1} \langle D_2H_d(\xk,\pko,\uk)-\Dt \xk,\delta
\pko \rangle \tau\nonumber\\
&&{+}\:\langle D_1H_d(\xk,\pko,\uk)+\Dt \pk,\delta \xk \rangle\tau  \langle+D_3H_d(\xk,\pko,\uk),\delta \uk \rangle\tau\nonumber\\
&&{+}\:\langle \phi_0 ,\delta
\lambda_0\rangle +\langle
\phi_n ,\delta \lambda_n\rangle +\langle -p_n+D_1\phi_n^T\lambda_n,\delta x_n\rangle +\langle
p_0+D_1\phi_0^T\lambda_0,\delta x_0\rangle\,,\nonumber\\
\end{eqnarray}
where the modified Leibnitz law (\ref{eq:modified_leibnitz}) has been used. We impose the variation
of the augmented cost function to be zero to obtain discrete necessary conditions for optimality
and transversality conditions:
\begin{eqnarray}
\Dt \xk&=&D_2H_d(\xk,\pko,\uk)\label{eq:dsneccond1}\,,\\
\Dt \pk&=&-D_1H_d(\xk,\pko,\uk)\label{eq:dsneccond2}\,,\\
0&=&D_3H_d(\xk,\pko,\uk)\label{eq:dsneccond3}\,,\\
p_0=-D_1\phi_0(x_0,t_0)^T\lambda_0\,,&\;& p_n=D_1\phi_n(x_n,t_n)^T\lambda_n\label{eq:dstranscond}\,.
\end{eqnarray}
The algorithm defined by (\ref{eq:dsneccond1}), (\ref{eq:dsneccond2}) and (\ref{eq:dsneccond3}) is
equivalent to the one derived by Bloch, Crouch, Marsden and Ratiu \cite{blo02cro} for the symmetric
rigid body.
\begin{lem}
The algorithm defined by (\ref{eq:dsneccond1}), (\ref{eq:dsneccond2}) and (\ref{eq:dsneccond3}) is
symplectic.
\end{lem}
\begin{pf}
Define the cost function $\bar{J}_a$ as:
\begin{equation}
\bar{J}_a=\sum_{k=0}^{n-1} (H_d(\xk,\pko,\uk)+\langle{\pko}, \Dt \xk\rangle)\tau\,.
\end{equation}
$\bar{J}_a$ is the augmented cost function from which we have removed the boundary conditions.
Boundary conditions yield transversality conditions, that is conditions on the initial and final
states of the system. Hence these terms are irrelevant to the study of the advance map
$(\xk,\pk,\uk)\mapsto (\xko,\pko,\uko)$. As in discrete dynamics, we consider $d^2 J_a$, assuming
$(\xk,\pk,\uk)$ verifies the above necessary conditions and we obtain:
\begin{equation}
d\bar{J}_a=\sum_{k=0}^{n-1}\Dt \langle \pk ,d \xk\rangle\tau\,.
\end{equation}
From $d^2=0$, we conclude:
\begin{equation}
0=\sum_{k=0}^{n-1}\Dt d\langle \pk ,d \xk\rangle\tau\,,\; \textrm{that is,}\;\forall k\in [0,n-1]
\,,\ d\pko\wedge d\xko=d\pk\wedge d\xk\,.
\end{equation}
\end{pf}
The symplectic nature of the algorithm is obtained directly from the variational principle - there
is no need to compute $d\pk \wedge d \xk$ and $d\pko \wedge d \xko$.

\paragraph{Midpoint scheme}
Midpoint discretization may also be obtained if we choose
$$
\begin{displaystyle}
\xdk=\frac{\xko+\xk}{2}\,,\,\pdk=\frac{\pko+\pk}{2}\,,\,\udk=\frac{\uko+\uk}{2}\,.
\end{displaystyle}
$$
and $\Dtd=R_{\tau/2}-R_{-\tau/2}$. One can readily verify that the discrete maximum principle
yields the following necessary conditions for optimality and transversality conditions:
\begin{eqnarray}
\Dtd \xdk&=&D_2H_d(\xdk,\pdk,\udk)\,,\label{eq:dneccond1_mp}\\
\Dtd \pdk&=&-D_1H_d(\xdk,\pdk,\udk)\,,\label{eq:dneccond2_mp}\\
0&=&D_3H_d(\xdk,\pdk,\udk)\,,\label{eq:dneccond3_mp}\\
p_0=D_1\phi_0(x_0,t_0)^T\lambda_0\,,&\;& p_n=-D_1\phi_n(x_n,t_n)^T\lambda_n\,.\label{eq:dtranscond_mp}
\end{eqnarray}
\begin{lem}
The algorithm defined by (\ref{eq:dneccond1_mp}), (\ref{eq:dneccond2_mp}) and
(\ref{eq:dneccond3_mp}) is symplectic.
\end{lem}
\begin{pf}
We omit the proof since it proceeds as before.
\end{pf}

\subsection{Discrete maximum principle v.s. discretization of the Pontryagin maximum principle}
So far we have considered two methods for obtaining a symplectic algorithm that integrates the necessary conditions for optimality. The first method, which applies only to a certain class of problems, consists of discretizing the necessary conditions obtained from the Pontryagin maximum principle once the control as been expressed as function of $(x,p$). The second method consists in using the discrete maximum principle. In this section, we show that under certain assumptions both methods are equivalent, that is we prove the commutative diagram (\ref{diag:1}).
\begin{equation}
\begin{diagram}
\node{\begin{array}{c}\min_u \int^T_0 g(x,u)dt\\\dot{x}=f(x,u)\end{array}}
\arrow{e}
\arrow{s,r}{(PMP)}
\node{\begin{array}{c}\min_u \sum_{k=0}^{n-1} g_d(\xdk,\udk)\\\Dtd \xdk=f_d(\xdk,\udk)\end{array}}
\arrow{s,r}{(DMP)}
\\
\node{\begin{array}{c}H(x,p,u)\\\bar{H}(x,p)\end{array}}
\arrow{e,t}{(DMHP)}
\node{\begin{array}{c}H_d(\xdk,\pdk,\udk)\\\bar{H}_d(\xdk,\pdk)\end{array}}
\label{diag:1}
\end{diagram}
\end{equation}
where $\bar{H}$ is defined by (\ref{eq:barh}), DMHP stands for discrete modified Hamilton's
principle, PMP stands for Pontryagin maximum principle, and DMP stands for discrete maximum
principle.

We recall the required assumptions to prove the equivalence of the diagram. We assume that (\ref{eq:neccond3}) can be
solved for $u$ as a function of $(x,p)$ and that the initial and final states $x(t_f)=x_f$ and
$x(t_0)=x_i$ are given. In addition, we impose $g_d=g$ and $f_d=f$.

To discretize the Hamiltonian system defined by $\bar{H}$, we use the discrete modified Hamilton's
principle:
\begin{equation}
0=\delta S^H_d=\delta \left(\tau\sum_{k=0}^{n-1} \langle\pdk, \Dtd\xdk\rangle
-\bar{H}(\xdk,\pdk)\right)
\label{eq:tmp4}
\end{equation}
for any variations of $(\xdk,\pdk)$ and $\delta x_0=\delta x_n=0$. One can readily check that
(\ref{eq:tmp4}) can also be written in an equivalent form as:
\begin{equation}
0=\delta S^H_d=\delta \left(\tau\sum_{k=0}^{n-1} \langle\pdk, \Dtd\xdk\rangle
-H(\xdk,\pdk,\udk)\right)
\end{equation}
for any variations of $(\xdk,\pdk,\udk)$ and $\delta x_0=\delta x_n=0$ where $\udk$ is now
considered as an independent variable. In addition since $f=f_d$ and $g=g_d$, $H=H_d$, and we
conclude that the discrete modified Hamilton's principle as formulated and the discrete maximum
principle are equivalent.

\subsection{The Heisenberg optimal control problem}
The Heisenberg problem (Brockett \cite{bro82}, Bloch et al. \cite{blo03bail}) refers to under actuated optimal control problems which are controllable. For instance, consider a particle that has two actuators in the $(x,y)$-plane and with velocity in the $z$ direction defined by $\dot{z}=y\dot{x}-x\dot{y}$. This system is controllable, however, to reach a point $(a>0,0,0)$ from the origin $(0,0,0)$ requires a non-trivial control vector. In the following, we study the Heisenberg problem to illustrate the approaches we have developed above. This problem formulates as:
\begin{equation}
\min_{u=(u_1,u_2)}\int_{t_0}^{t_f} \langle u,u \rangle dt\,,
\end{equation}
subject to
\begin{eqnarray}
\dot{x}&=&u\,,\\
\dot{y}&=&v\,,\\
\dot{z}&=&u y- v x\,,
\end{eqnarray}
and to the boundary conditions:
\begin{equation}
(x(t_0),y(t_0),z(t_0))=(0,0,0) \,, \  (x(t_f),y(t_f),z(t_f))=(a>0,0,0)\nonumber\,.
\end{equation}
This is a hard constraint problem, therefore the transversality conditions are of no use; They yield $2n$ equations but introduce $2n$ new variables.

Define $H$ as $$H(q,p,u)=\undemi \langle u,u\rangle + \langle p,\dot{q}\rangle\,,$$ where
$q=(x,y,z)$ and $p=(p_x,p_y,p_z)$. The Pontryagin maximum principle yields:
\begin{eqnarray}
\dot{q}&=&\pfrac{H}{p}(q,p,u)\label{eq:neccond1_heis}\,,\\
\dot{p}&=&-\pfrac{H}{q}(q,p,u)\label{eq:neccond2_heis}\,,\\
0&=&\pfrac{H}{u}(q,p,u)\,.
\label{eq:neccond3_heis}
\end{eqnarray}
Equation (\ref{eq:neccond3_heis}) allows us to solve for $u$ as a function of $(q,p)$:
\begin{eqnarray}
u_1=p_x+p_z y\,,\; u_2=p_y-p_z x\,,
\end{eqnarray}
Hence, equations (\ref{eq:neccond1_heis})-(\ref{eq:neccond2_heis}) become:
\begin{eqnarray}
\dot{q}&=&\pfrac{\bar{H}}{p}(q,p)\label{eq:neccond1bis_heis}\,,\\
\dot{p}&=&-\pfrac{\bar{H}}{q}(q,p)\label{eq:neccond2bis_heis}\,,
\end{eqnarray}
where
\begin{eqnarray}
\bar{H}(q,p)&=&H(q,p,u(q,p))\nonumber\\
&=&-\undemi (p_x^2+p_y^2)-p_x p_z y +p_y p_z x
\label{eq:barh_heis}
\end{eqnarray}
Equations (\ref{eq:neccond1bis_heis}) and (\ref{eq:neccond2bis_heis}) are of the same form as the Hamilton equations. Therefore, the necessary conditions for optimality yield a Hamiltonian system with Hamiltonian function $\bar{H}$. We now prove that $\bar{H}$ is degenerate at the origin, and so is the Legendre transform. The Hessian of $\bar{H}$ is:
\begin{equation}
\left(\pfrac{\bar{H}}{(q,p)}\right)=\begin{pmatrix} -1 & 0& -y \\
0&-1&x\\
-y&x&0
\end{pmatrix}\nonumber
\end{equation}
Thus, $\det\left(\pfrac{\bar{H}}{(q,p)}\right)=x^2+y^2$, i.e., the determinant of the Hessian of $\bar{H}$ is singular at $(0,0)$. As a result, it is not, \textit{a priori}, possible to define a Lagrangian function associated with the Hamiltonian $\bar{H}$ using the Legendre transform\footnote{ Using Lagrange multipliers one can define a Legendre transform and find a Lagrangian function associated with the system. We refer to Bloch \cite{blo03bail} for a presentation of this technique that involves variational principles with constraints.}. Therefore, the discrete modified Hamilton's principles (DMHP) must be used to discretize Eqns. \eqref{eq:neccond1bis_heis} and \eqref{eq:neccond2bis_heis}. One cannot use a discrete Hamilton's principles (DHP) for instance because the system is not Lagrangian. This point is of importance. It motivates the need to introduce the variational principles presented in this paper, as previous works on variational principles mostly focused on systems with non-degenerate Lagrangian functions.
To discretize the necessary conditions, we choose the geometry associated with the St{\"o}rmer rule and using DMHP (definition \ref{def:discrete_pcp_la}) to eventually find the following symplectic algorithm:
\begin{eqnarray}
\Dt \qk &=& D_2 \bar{H}(\qk,\pko)\,,\\
\Dt \pk &=&-D_1 \bar{H}(\qk,\pko)\,.
\end{eqnarray}

Let us now discretize the Heisenberg problem using the second approach, based on the use of the discrete maximum principle. We first discretize the problem statement:
\begin{equation}
\min_{\uk=(\uok,\utk)} \undemi\sum_{k=0}^{n-1} \langle \uk,\uk \rangle\,,
\end{equation}
subject to
\begin{eqnarray}
\Dt \xk&=&\uok \,,\\
\Dt \yk&=&\utk \,,\\
\Dt \zk&=&\uok \yk -\utk \xk\,.
\end{eqnarray}
Define  the discrete augmented cost function $J_a$:
\begin{equation}
J_a=\sum_{k=0}^{n-1} H_d(\qk,\pko,\uk)-\langle \pko, \Dt \qk\rangle\,,
\end{equation}
where $H_d(\qk,\pko,\uk)=\langle \uk,\uk \rangle+\langle \pko,\qk \rangle$ and $\qk=(\xk,\yk,\zk)$.
To find discrete necessary conditions for optimality we set the variations of $J_a$ to zero, and we
obtain:
\begin{eqnarray}
\Dt \qk&=&D_2H_d(\qk,\pko,\uk)\label{eq:neccond1_dheis}\,,\\
\Dt \pk&=&-D_1H_d(\qk,\pko,\uk)\label{eq:neccond2_dheis}\,,\\
0&=&D_3H_d(\qk,\pko,\uk)\label{eq:neccond3_dheis}\,.
\end{eqnarray}
Equation (\ref{eq:neccond1_dheis}) allows us to find $\uk$ as a function of $(\qk,\pko)$:
\begin{equation}
\uok=p_{x,k+1} + p_{z,k+1} \yk\,,\;\utk=p_{y,k+1} - p_{z,k+1} \xk\,.
\end{equation}
We then substitute these expressions into equations
(\ref{eq:neccond1_dheis})-(\ref{eq:neccond2_dheis}):
\begin{eqnarray}
\Dt \qk&=&D_2\bar{H}_d(\qk,\pko)\label{eq:neccond1bis_dheis}\,,\\
\Dt \pk&=&-D_1\bar{H}_d(\qk,\pko)\label{eq:neccond2bis_dheis}\,,
\end{eqnarray}
where $\bar{H}_d(\qk,\pko)=H_d(\qk,\pko,\uk(\qk,\pko))$. By virtue of the commutative diagram,
(\ref{eq:neccond1bis_dheis}) and (\ref{eq:neccond2bis_dheis}) define the same symplectic algorithm
as (\ref{eq:neccond1bis_heis}) and (\ref{eq:neccond2bis_heis}).
\newline

In this example, we chose a trivial discretization of the dynamics and of the cost function; $f=f_d$
and $g=g_d$. Other algorithms may be obtained using  nontrivial discretizations. In that case the
equivalence principle may not hold but the algorithm we obtain will still be symplectic. In
addition, in this example we did not take into account any boundary conditions since we have seen
earlier in the paper that both methods yield comparable transversality conditions. Finally, as in
discrete dynamics, the discrete maximum principle may be modified in order to yield
symplectic-energy conserving algorithms. We add an independent parameter $\tau_k$ and consider the
time as a generalized coordinate, the optimal control problem then formulates as follows:
\begin{equation}
min_u \sum_{k=0}^{n-1} g_d(\xdk,\udk) (\tko-\tk) =\sum_{k=0}^{n-1} g_d(\xdk,\udk) \Dt\tk \tau\,.
\end{equation}
subject to the dynamics
\begin{equation}
\Dtd \xdk=\Dtd \tdk f_d(\xdk,\udk)\,.
\label{eq:dgdyn}
\end{equation}

\section{Conclusions}
In this paper we have presented a general framework that allows one to study discrete systems. We
have introduced variational principles on the tangent and cotangent bundles  that are the discrete
counterpart of the known principles of critical action for Lagrangian and Hamiltonian dynamical
systems. We have shown that they allowed us to recover most of the classical symplectic algorithms.
In the future, we will try to derive additional symplectic algorithms such as the symplectic
partitioned Runge-Kutta algorithm. In addition, we have seen that by increasing the dimensionality
of the configuration space, symplectic algorithms may be transformed into symplectic-energy
conserving algorithms. When time is a generalized coordinate, the dynamical system is subject to an
energy constraint and we are able to adapt our variational principles to take into account such a
constraint. In the same manner, our approach may be modified to derive symplectic algorithms to
integrate non-autonomous dynamical and control systems with (non-holonomic) constraints. We have
also identified a class of coordinate transformations that leaves the variational principles
presented in this paper invariant and developed a discrete Hamilton-Jacobi theory. This theory
allows us to relate the energy error in the integration using different set of coordinates.
Finally, for optimal control problems we have developed a discrete maximum principle that yields
discrete necessary conditions for optimality. These conditions are in agreement with the usual
conditions obtained from Pontryagin maximum principle. In future research, we want to use the
general framework introduced in this paper to develop variational principles for multi-symplectic
algorithms, that is a spacetime discretization will be used instead of the time discretization.
Such a formulation would allows us to develop efficient numerical algorithms for simulation of the
motion of rigid bodies and complex interconnected systems.

{\bf Acknowledgement:} We would like to thank Jerry Marsden for valuable discussions.

\bibliographystyle{plain}
\bibliography{control,discrete,dynamic}

\end{document}